\documentclass[graybox,envcountsect]{svmult}

\usepackage{graphicx}
\usepackage{amssymb}
\usepackage{amsmath}
\usepackage{color}

\usepackage[latin1]{inputenc}
\usepackage{amsfonts}
\usepackage{tikz} 
\usetikzlibrary{matrix}
\usepackage[english]{babel}
\usepackage[latin1]{inputenc}
\usepackage[T1]{fontenc}
\usepackage{ae}

\usepackage{url}
\usepackage{hyperref}

\usepackage{color}
\usepackage{tikz}
\usetikzlibrary{matrix}
\usepackage{cite}           
\usepackage{mathptmx}       
\usepackage{helvet}         
\usepackage{courier}        
\usepackage{type1cm}        
\usepackage{makeidx}         
\usepackage{graphicx}        
\usepackage{multicol}        
\usepackage[bottom]{footmisc}

\newtheorem{cor}[theorem]{Corollary}
\definecolor{shadecolor}{rgb}{0.8, 0.8, 0.8}
\newcommand\eqn[1]{(\ref{eq:#1})}
\def\maple{\textsc{maple \ }}
\def\maplep{\textsc{maple}.\ }

\def\carrot{${}^\wedge$}
\newcommand\exampleshade[1]{
\begin{shaded}
\null\noindent
EXAMPLE:\\
{\tt
\noindent
#1
}
\end{shaded}
}
\newcommand\mapleshade[1]{
\begin{shaded}
\null\noindent
\hphantom{EXAMPLE:}\\
{\tt
\noindent
#1
}
\end{shaded}
}
\newcommand{\SL}{\mbox{SL}}
\newcommand{\SLZ}{\mbox{SL}_2(\mathbb{Z})}

\newcommand\mylabel[1]{\label{#1}}
\newcommand{\beqs}{\begin{equation*}}
\newcommand{\eeqs}{\end{equation*}}
\newcommand{\beq}{\begin{equation}}
\newcommand{\eeq}{\end{equation}}
\newcommand\thm[1]{\ref{thm:#1}}

\newcommand\corol[1]{\ref{cor:#1}}
\newcommand\sect[1]{\ref{sec:#1}}
\newcommand\subsect[1]{\ref{subsec:#1}}

\newcommand{\C}{\mathbb{C}}

\newcommand{\Z}{\mathbb{Z}}
\newcommand{\Hup}{\mathcal{H}}

\renewcommand{\Im}{\mathrm{Im}\,}
\newcommand{\leg}[2]{\genfrac{(}{)}{}{}{#1}{#2}}

\newcommand\ord{\mbox{ord}\,}
\newcommand\ORD{\mbox{ORD}\,}

\makeindex             

\begin{document}
\title*{Automatic Proof of Theta-Function Identities}

\author{Jie Frye and Frank Garvan 
\thanks{ A preliminary version of this paper was presented by J.\ Frye
on January 10, 2013 at JMM2013, San Diego. 
F.\ Garvan was supported in part by a grant from
the Simon's Foundation (\#318714). }}
\institute{Jie Frye,\\
Bunker Hill Community College,	Boston, MA 02129\\
\email{jiefrye@gmail.com}\\
Frank Garvan,\\
Department of Mathematics, University of Florida, Gainesville, FL 32601\\
\email{fgarvan@ufl.edu}}

\maketitle
\date{July 20, 2018}

\abstract{
This is a tutorial for using two new MAPLE packages, \texttt{thetaids}
and \texttt{ramarobinsids}. The \texttt{thetaids} package is designed for proving
generalized eta-product identities using the valence formula
for modular functions. We show how this package can be
used to find theta-function identities as well as prove them. As
an application, we show how to find and prove Ramanujan's 40
identities for his so called Rogers-Ramanujan functions $G(q)$ and
$H(q)$. In his thesis Robins found similar identities for higher
level generalized eta-products. Our \texttt{ramarobinsids} package is for
finding and proving identities for generalizations of Ramanujan's
$G(q)$ and $H(q)$ and Robin's extensions. These generalizations are
associated with certain real Dirichlet characters. We find a total
of over 150 identities.
}


%
%
\renewcommand{\theequation}{\arabic{section}.\arabic{equation}}
\allowdisplaybreaks

\section{Introduction}
\setcounter{equation}{0}
\mylabel{sec:intro}

The Rogers-Ramanujan functions are
\begin{align}
G(q) &= \sum_{n=0}^\infty \frac{q^{n^2}}{(q; q)_n} = 
      \prod_{n=0}^\infty \frac{1}{(1-q^{5n+1})(1-q^{5n+4})},
\mylabel{eq:RRfunc}\\
H(q) &= \sum_{n=0}^\infty \frac{q^{n(n+1)}}{(q; q)_n} = 
      \prod_{n=0}^\infty \frac{1}{(1-q^{5n+2})(1-q^{5n+3})}.
\nonumber
\end{align}
The ratio of these two functions is the famous Rogers-Ramanujan continued
fraction \cite{An81c}
\begin{align*}
\frac{G(q)}{H(q)} &= 
\prod _{n=0}^{\infty }{\frac {(1-q^{5n+2})(1-q^{5n+3})}{(1-q^{5n+1})(1-q^{5n+4})}}\\
&=1 + {\cfrac
{q}{1+{\cfrac {q^2}{1+{\cfrac {q^{3}}{1+{\cfrac
{q^{4}}{1+\ddots }}}}}}}}
\end{align*}               

Ramanujan also found
\begin{equation} \mylabel{eq:ENTRY3.1}
H(q)G(q)^{11} - q^2G(q)H(q)^{11} = 1+11G(q)^6H(q)^6
\end{equation}
and
\begin{equation} \mylabel{eq:ENTRY3.4}
H(q)G(q^{11}) - q^2G(q)H(q^{11}) = 1,
\end{equation}
and remarked that ``each of these formulae is the simplest of a large class". 
Here we have used the standard $q$-notation
\[ (a; q)_n := \prod_{j=0}^{n-1} (1-aq^j) \qquad  
(a; q)_\infty := \prod_{j=0}^{\infty} (1-aq^j). 
\]
In 1974 B.~J.~Birch \cite{Bi74} published a description of some manuscripts of 
Ramanujan
including a list of forty identities for the Rogers-Ramanujan functions.
Biagioli \cite{Bi89} show how the theory of modular forms could prove
identities of this type efficiently. See \cite{Be-Ye} and \cite{Be-etal07}
for recent work. It is instructive to write the Rogers-Ramanujan functions
in terms of generalized eta-products.

The Dedekind eta-function is defined by
\[ 
\eta(\tau) = q^{\frac{1}{24}} \prod_{n=1}^\infty (1-q^n),
\]
where $\tau \in \Hup:= \{ \tau \in \C : \Im \tau > 0\}$ and $q := e^{2\pi i \tau}$, 
and the generalized Dedekind eta function is defined to be 
\begin{equation}
\eta_{\delta;g}(\tau) =q^{\frac{\delta}{2} P_2(g/\delta) }
\prod_{m\equiv \pm g\pmod{\delta}} (1 - q^m),
\mylabel{eq:Geta}
\end{equation}
where $P_2(t) = \{t\}^2 - \{t\} + \tfrac16$ is the second periodic Bernoulli polynomial, 
$\{t\}=t-[t]$ is the fractional part of $t$, $g,\delta,m\in \mathbb{Z}^{+}$
and $0 \le g < \delta$. The function $\eta_{\delta;g}(\tau)$ is a modular function
(modular form of weight $0$) on $\SLZ$ with a multiplier system.

Ramanujan's identities \eqn{ENTRY3.1} and \eqn{ENTRY3.4} can be rewritten as
\beq
\frac{1}{\eta_{5;2}(\tau) \eta_{5;1}(\tau)^{11}}
-
\frac{1}{\eta_{5;1}(\tau) \eta_{5;2}(\tau)^{11}}  = 1 + 11\,
\frac{\eta(5\tau)^{6}}{\eta(\tau)^{6}}
\mylabel{eq:altENTRY3.1}
\eeq
and
\beq
\frac{1}{\eta_{5;2}(\tau) \eta_{5;1}(11\tau)}
-
\frac{1}{\eta_{5;1}(\tau) \eta_{5;2}(11\tau)}  = 1.
\mylabel{eq:altENTRY3.4}
\eeq
It is natural to consider higher level analogues of Ramanujan's identities
\eqn{ENTRY3.1} and \eqn{ENTRY3.4}. The following are nice level $13$ analogues:
\beq
\frac{1}{\eta_{13;2,5,6}(\tau) \eta_{13;1,3,4}(\tau)^{3}}
-
\frac{1}{\eta_{13;1,3,4}(\tau) \eta_{13;2,5,6}(\tau)^{3}}  = 1 + 3\,
\frac{\eta(13\tau)^{2}}{\eta(\tau)^{2}}
\mylabel{eq:RR13a}
\eeq
and
\beq
\frac{1}{\eta_{13;2,5,6}(\tau) \eta_{13;1,3,4}(3\tau)}
-
\frac{1}{\eta_{13;1,3,4}(\tau) \eta_{13;2,5,6}(3\tau)}  = 1.          
\mylabel{eq:RR13b}
\eeq
Here we have used the notation
$$
\eta_{\delta;g_1,g_2,\dots,g_k}(\tau) =
\eta_{\delta;g_1}(\tau) \,
\eta_{\delta;g_1}(\tau) \,
\cdots
\eta_{\delta;g_k}(\tau).
$$
Equation \eqn{RR13a} was found by Ramanujan \cite[Eq.(8.4),p.373]{BeRNIII}, and equation \eqn{RR13b}
is due to Robins \cite{Ro-thesis}, who considered more general identities.
The following is level $17$ analogue of \eqn{RR13b} and appears to be new.
\beq
\frac{1}
{\eta_{17;3,5,6,7}(\tau) \eta_{17;1,2,4,8}(2\tau)}
-
\frac{1}
{\eta_{17;1,2,4,8}(\tau) \eta_{17;3,5,6,7}(2\tau)} = 1.
\mylabel{eq:RR17b}
\eeq
Motivated by these examples and other work of Robins \cite{Ro-thesis}
one is led naturally to consider
\beq
G(n,N,\chi) = G(n) := 
\prod_{\substack{\chi(g) = 1\\ 0 < g < \frac{N}{2}}}
\frac{1}{\eta_{N;g}(n\tau)}, \quad
H(n,N,\chi) = H(n) := 
\prod_{\substack{\chi(g) = -1\\ 0 < g < \frac{N}{2}}}
\frac{1}{\eta_{N;g}(n\tau)},
\mylabel{eq:GHdefs}
\eeq
where $\chi$ is a non-principal real Dirichlet character mod $N$ satisfying
$\chi(-1)=1$. Ratios of functions of this type were studied by
Huber and Schultz \cite{Hu-Sc16}. They found the following level $17$
identity:
\beq
(r^2+8\,r-1)^2\,s^2-2\,r\,(r^2+1)\,s+r^2 = 0,
\mylabel{eq:HuSc17id}
\eeq
where 
$$
r = \frac{ H(1,17,\leg{\cdot}{17})}
         { G(1,17,\leg{\cdot}{17})},\quad
s = \frac{\eta(17\tau)^3}{\eta(\tau)^3}.
$$

The main goal of the \texttt{thetaids} \maple package is
to 
automatically prove identities for generalized eta-products 
using the  theory of modular functions.
In Sections \sect{genramarobins}-\sect{moregenramarobins} 
we describe the \texttt{ramarobinsids}
package, which uses the \texttt{thetaids} package to search for and prove
identities for general functions $G(n,N,\chi)$ and $H(n,N,\chi)$ that
are like the 
theta-function identities considered by Ramanujan \cite{Be-etal07}
and Robins \cite{Ro-thesis}. 

We note that Liangjie \cite{Ye17} gave an algorithm for proving
relations for certain theta-functions and their derivatives using a 
different method.
We also note that Lovejoy and Osburn 
\cite{Lo-Os13a},  \cite{Lo-Os13b},  \cite{Lo-Os13c},  \cite{Lo-Os15},  
have used an earlier version
of the \texttt{thetaids} package to prove theta-functions identities that were 
needed to 
establish an number of results for mock-theta functions.

\subsection{Installation Instructions}

    First install the \texttt{qseries} package from 
\begin{center}
\url{http://qseries.org/fgarvan/qmaple/qseries}
\end{center}
 and  follow the directions on that page. Before proceeding it is advisable
to become familiar with the functions in the \texttt{qseries} package.
See \cite{Ga99b} for a tutorial.
Then go to
\begin{center}
\url{http://qseries.org/fgarvan/qmaple/thetaids}
\end{center}
 to install the \texttt{thetaids}
package. In Section \sect{genramarobins} you will need to install
the \texttt{ramarobinsids} package from
\begin{center}
\url{http://qseries.org/fgarvan/qmaple/ramarobinsids}
\end{center}

\section{Proving theta-function identities}
\mylabel{sec:proving}
\setcounter{equation}{0}

          To prove a given theta-function identity one needs to  basically
do the following.
\begin{enumerate}
\item[(i)]
Rewrite the identity in terms of generalized eta-functions.
\item[(ii)]
Check that each term in the identity is a modular function on some
group $\Gamma_1(N)$.
\item[(iii)]
Determine the order at each cusp of $\Gamma_1(N)$ of each term
in the identity.
\item[(iv)]
Use the valence formula to determine up to which power of $q$ is needed to verify
the identity.
\item[(v)]
Finally prove the identity by carrying out the verification.
\end{enumerate}

In this section we explain how to carry out each of these steps in \maplep
Then we show how the whole process of proof can be automated.

\subsection{Encoding theta-functions, eta-functions and generalized 
eta-functions}
\mylabel{subsec:encoding}

We recall Jacobi's triple product for theta-functions:
\beq
\prod_{n=1}^\infty (1 - zq^{n-1}) (1 - z^{-1}q^n) (1 - q^n)
=
\sum_{n=-\infty}^\infty (-1)^n z^n q^{n(n-1)/2},
\mylabel{eq:JTP}
\eeq
so that
\beq
\prod_{n=1}^\infty(1-q^{\delta n-g})(1-q^{\delta n+g-\delta})
(1 - q^{\delta n})
=
\sum_{n=-\infty}^\infty (-1)^n 
q^{\tfrac{1}{2}n(\delta n-\delta+2g)}.
\mylabel{eq:JTP2}
\eeq
In the \texttt{qseries} \maple package the function on the left side of
\eqn{JTP2} is encoded symbolically as 
\texttt{JAC(g,$\delta$,infinity)}. This is the building block of the functions
in our package. In the \texttt{qseries} package \texttt{JAC(0,$\delta$,infinity)}
corresponds symbolically to 
\beq
\prod_{n=1}^\infty(1-q^{\delta n})
=
\sum_{n=-\infty}^\infty (-1)^n 
q^{\tfrac{\delta}{2}n(3 n + 1)},
\mylabel{eq:EPNT}
\eeq
which is Euler's Pentagonal Number Theorem.

\begin{center}
\begin{tabular}{|c|c|}
\noalign{\hrule}
Function &  Symbolic \maple form \\
\noalign{\hrule}
$\displaystyle
\prod_{n=1}^\infty(1-q^{\delta n-g})(1-q^{\delta n+g-\delta})
(1 - q^{\delta n})$
& \texttt{JAC(g, $\delta$, infinity)} \\
\noalign{\hrule}
$\displaystyle\prod_{n=1}^\infty(1-q^{\delta n})$
 & \texttt{JAC(0, $\delta$, infinity)} \\
\noalign{\hrule}
$\eta_{\delta;g}(\tau)$ & \texttt{GETA($\delta$,g)} \\         
\noalign{\hrule}
$\eta(\delta\tau)$ & \texttt{EETA($\delta$)} \\         
\noalign{\hrule}
\end{tabular}
\end{center}

We will also consider generalized eta-products. Let $N$ be a fixed positive integer.
A generalized Dedekind eta-product of level $N$ has the form
\begin{equation} \mylabel{eq:fdef}
f(\tau) = \prod_{\substack{\delta\mid N\\ 0 < g < \delta}}
                 \eta_{\delta;g}^{r_{\delta,g}}(\tau),
\end{equation}
where
\begin{equation} \mylabel{r-dg}
     r_{\delta,g} \in
     \begin{cases}
      \frac{1}{2}\Z & \mbox{if $g=\delta/2$},\\
      \Z & \mbox{otherwise}.
      \end{cases}
\end{equation}

In \maple we represent the generalized eta-product
$$
\eta_{N;g_1}(\tau)^{r_1} \, \eta_{N;g_2}(\tau)^{r_2} \, \cdots
\eta_{N;g_m}(\tau)^{r_m} 
$$
symbolically by the list
$$
[[N,g_1,r_1], [N,g_2,r_2]\, \dots, [N,g_m,r_m]].
$$
We call such a list a \textit{geta-list}.

\subsection{Symbolic product conversion} 
\mylabel{subsec:symprodcon}

\hphantom{X}

\texttt{jac2eprod} --- Converts a quotient of theta-functions in \texttt{JAC}
notation to a product of generalized eta-functions in \texttt{EETA} and \texttt{GETA}
notation.

\exampleshade{
> with(qseries):\\
> with(thetaids):\\
> G:=q->add(q\carrot(n\carrot2)/aqprod(q,q,n),n=0..10):\\
> H:=q->add(q\carrot(n\carrot2+n)/aqprod(q,q,n),n=0..10):\\
> JG:=jacprodmake(G(q),q,50);
$$
                    \frac{JAC(0, 5, \infty)}{JAC(1, 5, \infty)}
$$
> HG:=jacprodmake(H(q),q,50);
$$
                      \frac{JAC(0, 5, \infty)}{JAC(2, 5, \infty)}
$$
> JP:=jacprodmake(H(q)*G(q)\carrot(11),q,80);
$$
{\frac { \left( {\it JAC} \left( 0,5,\infty  \right)  \right) ^{12}}{
 \left( {\it JAC} \left( 1,5,\infty  \right)  \right) ^{11}{\it JAC}
 \left( 2,5,\infty  \right) }}
$$
> GP:=jac2eprod(JP);
$$
{\frac {1}{ \left( {\it GETA} \left( 5,1 \right)  \right) ^{11}{\it 
GETA} \left( 5,2 \right) }}
$$
}

\texttt{jac2getaprod} --- Converts a quotient of theta-function in \texttt{JAC}
notation to a product of generalized eta-functions in standard notation.

\exampleshade{
> jac2getaprod(JP);
$$
\frac {1}{  \eta_{5,1}  \left( \tau \right)^{11} 
\eta_{5,2}\left( \tau \right)}
$$
}

\texttt{GETAP2getalist} --- Converts 
a product of generalized eta-functions into a list as described above.

\exampleshade{
> GETAP2getalist(GP);
$$
                   [[5, 1, -11], [5, 2, -1]]
$$
}

\subsection{Processing theta-functions} 
\mylabel{subsec:processing}

There are two main functions in the \texttt{thetaids} package for
processing combinations of theta-functions.

\texttt{mixedjac2jac} --- Converts a sum of quotients of theta-functions
written in terms of \texttt{JAC(a,b,infinity)} to a sum with
the same base \texttt{b}. The functions \texttt{jac2series} and \texttt{jacprodmake}
from the \texttt{qseries} package are used.

\exampleshade{
> Y1:=1+jacprodmake(G(q),q,100)*jacprodmake(H(q\carrot2),q,100);
$$
1+{\frac {{\it JAC} \left( 0,5,\infty  \right) {\it JAC} \left( 0,10,
\infty  \right) }{{\it JAC} \left( 1,5,\infty  \right) {\it JAC}
 \left( 4,10,\infty  \right) }}
$$
> Y2:=mixedjac2jac(Y1);
$$
1+{\frac { \left( {\it JAC} \left( 0,10,\infty  \right)  \right) ^{3}
}{{\it JAC} \left( 1,10,\infty  \right)  \left( {\it JAC} \left( 4,10,
\infty  \right)  \right) ^{2}}}
$$
}

\texttt{processjacid} --- Processes a theta-function identity
written as a rational function of \texttt{JAC}-functions 
using \texttt{mixedjac2jac} and
renormalizing by dividing by the term with the lowest power of $q$.

As an example, we consider the well-known identity
\beq
\theta_3(q)^4 = \theta_2(q)^4 + \theta_4(q)^4.
\mylabel{eq:thetaid1a}
\eeq

\exampleshade{
> with(qseries): \\
> with(thetaids): \\
> F1:=theta2(q,100)\carrot4: \\
> F2:=theta3(q,100)\carrot4: \\
> F3:=theta4(q,100)\carrot4: \\
> findhom([F1,F2,F3],q,1,0);
$$
 \left\{ X_{{1}}-X_{{2}}+X_{{3}} \right\} 
$$
> JACID0:=qs2jaccombo(F1-F2+F3,q,100);
$$
16\,{\frac {q \left( {\it JAC} \left( 0,4,\infty  \right)  \right) ^{6
}}{ \left( {\it JAC} \left( 2,4,\infty  \right)  \right) ^{2}}}-{
\frac { \left( {\it JAC} \left( 0,4,\infty  \right)  \right) ^{6}
 \left( {\it JAC} \left( 2,4,\infty  \right)  \right) ^{6}}{ \left( {
\it JAC} \left( 1,4,\infty  \right)  \right) ^{8}}}+ \left( {\it JAC}
 \left( 1,2,\infty  \right)  \right) ^{4}
$$
> JACID1:=processjacid(JACID0);
$$
-16\,{\frac {q \left( {\it JAC} \left( 1,4,\infty  \right)  \right) ^{
8}}{ \left( {\it JAC} \left( 2,4,\infty  \right)  \right) ^{8}}}+1-{
\frac { \left( {\it JAC} \left( 1,4,\infty  \right)  \right) ^{16}}{
 \left( {\it JAC} \left( 0,4,\infty  \right)  \right) ^{12} \left( {
\it JAC} \left( 2,4,\infty  \right)  \right) ^{4}}}
$$
> expand(jac2getaprod(JACID1));
$$
-{\frac{\eta_{4;1}(\tau)^{16}}{\eta_{4;2}(\tau)^{4}}}
+1
-16\,{\frac{\eta_{4;1}(\tau)^{8}}{\eta_{4;2}(\tau)^{8}}}
$$
}

We see that \eqn{thetaid1a} is equivalent to the identity
\beq
{\frac{\eta_{4;1}(\tau)^{16}}{\eta_{4;2}(\tau)^{4}}}
+
16\,{\frac{\eta_{4;1}(\tau)^{8}}{\eta_{4;2}(\tau)^{8}}}
=1.
\mylabel{eq:thetaid1b}
\eeq

\subsection{Checking modularity} 
\mylabel{subsec:modcheck}

Robins \cite{Ro94} has found sufficient conditions under which a 
generalized eta-product is a modular function on $\Gamma_1(N)$.
\begin{theorem}[\cite{Ro94}(Theorem 3)]
\mylabel{thm:gepmodfunc}
The function $f(\tau)$, defined in \eqn{fdef}, is a 
modular function on $\Gamma_1(N)$ if
\begin{enumerate}
\item[(i)]
$\displaystyle\sum_{\substack{\delta \mid N \\ g}}
\delta P_2(\textstyle{\frac{g}{\delta}}) r_{\delta,g} \equiv 0 \pmod{2}$, 
and
\item[(ii)]
$\displaystyle\sum_{\substack{\delta \mid N \\ g} }
\frac{N}{\delta} P_2(0) r_{\delta,g} \equiv 0 \pmod{2}$.
\end{enumerate}
\end{theorem}

The functions on the left side of (i), (ii) above are computed
using the \maple functions \texttt{vinf} and \texttt{v0} respectively.
Suppose $f(\tau)$ is given as in \eqn{fdef} and this generalized
eta-product is encoded as the geta-list $L$. Recall that each item in the
list $L$ has the form $[\delta,g,r_{\delta,g}]$. The syntax is
\texttt{vinf(L,N)} and \texttt{v0(L,N)}. As an example we consider the two
generalized eta-products in \eqn{thetaid1b}.

\exampleshade{
> L1:=[[4,1,16],[4,2,-4]];
$$
                    [[4, 1, 16], [4, 2, -4]]
$$
> vinf(L1,4),v0(L1,4);
$$
                              0, 2
$$
> L2:=[[4,1,8],[4,2,-8]];
$$
                    [[4, 1, 8], [4, 2, -8]]
$$
> vinf(L2,4),v0(L2,4);
$$                              
2, 0
$$
}
The numbers $0$, $2$ are even and we see that both generalized eta-products
in \eqn{thetaid1b} are modular functions on $\Gamma_1(4)$
by Theorem \thm{gepmodfunc}.

\texttt{Gamma1ModFunc(L,N)} --- Checks whether a given generalized
eta-product is a modular function on $\Gamma_1(N)$. Here the generalized
eta-product is encoded as the geta-list $L$. The function first checks
whether each $\delta$ is a divisor of $N$ and checks whether both
\texttt{vinf(L,N)} and \texttt{v0(L,N)} are even. It returns $1$ if it is
a modular function on $\Gamma_1(N)$ otherwise it returns $0$. If the global
variable \texttt{xprint} is set to \textit{true} then more detailed
information is printed. Thus here and throughout \texttt{xprint} can be
used for debugging purposes.

\exampleshade{
> Gamma1ModFunc(L1,4);
$$
1
$$
> xprint := true: \\
> Gamma1ModFunc(L1,4); \\
* starting Gamma1ModFunc with L=[[4, 1, 16], [4, 2, -4]] and N=4 \\
All n are divisors of 4\\
val0=2\\
which is even.\\
valinf=0\\
which is even.\\
It IS a modfunc on Gamma1(4)
$$
                               1
$$
}

\subsection{Cusps}
\mylabel{subsec:cusps}

Cho, Koo and Park \cite{Ch-Ko-Pa} have found a set of inequivalent 
cusps for $\Gamma_1(N) \cap \Gamma_0(mN)$. 
The group $\Gamma_1(N)$ corresponds to the case $m=1$.
\begin{theorem}[\cite{Ch-Ko-Pa}(Corollary 4, p.930)]
\mylabel{thm:chokoopark}
Let $a$, $c$, $a'$, $c\in\Z$ with $(a,c)=(a',c')=1$.
\begin{enumerate}
\item[(i)]
The cusps $\frac{a}{c}$ and $\frac{a'}{c'}$ are equivalent mod $\Gamma_1(N)$ 
if and only if
$$
\begin{pmatrix}
a' \\ c'
\end{pmatrix}
\equiv \pm
\begin{pmatrix}
 a + nc\\
 c
\end{pmatrix}
\pmod{N}
$$
for some integer $n$.
\item[(ii)]
The following is a complete set of inequivalent cusps mod $\Gamma_1(N)$.
\begin{align*}
\mathcal{S} &= \left\{ \frac{y_{c,j}}{x_{c,i}} \,:\, 0 < c \mid N,\,
  0 < s_{c,i},\, a_{c,j} \le N,\, (s_{c,i},N)=(a_{c,j},N)=1,\right.\\
     &\qquad s_{c,i}=s_{c,i'} \iff s_{c,1}\equiv\pm s_{c',i'} 
                  \pmod{{\textstyle \frac{N}{c}}},\\
    &\qquad a_{c,j}=a_{c,j'} \iff
\begin{cases}
a_{c,j}\equiv\pm a_{c,j'} \pmod{c}, &\mbox{if $c =\frac{N}{2}$ or $N$},\\
a_{c,j}\equiv a_{c,j'} \pmod{c}, &\mbox{otherwise},
\end{cases} \\
&\left. \vphantom{ \frac{y_{c,j}}{x_{c,i}} } x_{c,i}, y_{c,j}\in\Z\,
\mbox{chosen s.th.}\,
x_{c,i}\equiv c s_{c,i},\, y_{c,j}\equiv a_{c,j}\pmod{N},
\, (x_{c,i},y_{c,j})=1\right\},
\end{align*}
\item[(iii)]
and the fan width of the cusp $\frac{a}{c}$ is given by
$$
\kappa({\textstyle \frac{a}{c}},\Gamma_1(N)) =
\begin{cases}
1, & \mbox{if $N=4$ and $(c,4)=2$},\\
\frac{N}{(c,N)}, & \mbox{otherwise}.
\end{cases}
$$
\end{enumerate}
\end{theorem}
In this theorem, it is understood as usual that the fraction $\frac{\pm1}{0}$ 
corresponds to $\infty$.

\texttt{cuspequiv1}$(a_1,c_1,a_2,c_2,N)$ --- determines whether the cusps $a_1/c_1$
and $a_2/c_2$ are $\Gamma_1(N)$-equivalent using Theorem \thm{chokoopark}(i).

\exampleshade{
> cuspequiv1(1,3,1,9,40);
$$
                             false
$$
> cuspequiv1(1,9,2,9,40);
$$
                              true
$$
}
We see that modulo $\Gamma_1(40)$ the cusps $\tfrac13$ and $\tfrac19$ are inequivalent
and the cusps $\tfrac19$ and $\tfrac29$ are equivalent.

\texttt{Acmake(c,N)} ---  returns the set $\{a_{c,j}\}$ where $c$ is a positive
divisor of $N$.

\texttt{Scmake(c,N)} ---  returns the set $\{s_{c,i}\}$ where $c$ is a positive
divisor of $N$.

\texttt{newxy(x,y,N)} --- returns $[x_1,y_1]$ for given $(x,y,N)=1$ such that
$x_1\equiv x \pmod{N}$ and
$y_1\equiv y \pmod{N}$.             

\texttt{cuspmake1(N)} --- 
   returns a set of inequivalent cusps                   
   for $\Gamma_1(N)$ using Theorem \thm{chokoopark}.
   Each cusp $a/c$ in the list is represented by $[a,c]$,
   so that $\infty$ is represented by $[1,0]$. This \maple procedure uses the
   functions \texttt{Acmake}, \texttt{Scmake} and \texttt{newxy}.

\texttt{cuspwid1(a,c,N)} ---                                     
   returns the width of the cusp $a/c$
   for the group $\Gamma_1(N)$ using Theorem \thm{chokoopark}(iii).

\exampleshade{
> C10:=cuspmake1(10);
$$
\{[0, 1], [1, 0], [1, 2], [1, 3], [1, 4], [1, 5], [2, 5], [3, 10]\}
$$
> for L in C10 do lprint(L,cuspwid1(L[1],L[2],10));od;\\
\hphantom{POO} [0, 1], 10 \\
\hphantom{POO} [1, 0], 1\\
\hphantom{POO} [1, 2], 5\\
\hphantom{POO} [1, 3], 10\\
\hphantom{POO} [1, 4], 5\\
\hphantom{POO} [1, 5], 2\\
\hphantom{POO} [2, 5], 2\\
\hphantom{POO} [3, 10], 1\\
}

We have the following table of cusps for $\Gamma_1(10)$.
\begin{center}
\begin{tabular}{|c|c|}
\noalign{\hrule}
cusp & cusp-width\\
\noalign{\hrule}
$0$  & $10$ \\
\noalign{\hrule}
$\infty$  & $1$ \\
\noalign{\hrule}
$\tfrac12$ & $5$ \\
\noalign{\hrule}
$\tfrac13$ & $10$ \\
\noalign{\hrule}
$\tfrac14$ & $5$ \\
\noalign{\hrule}
$\tfrac15$ & $2$ \\
\noalign{\hrule}
$\tfrac25$ & $2$ \\
\noalign{\hrule}
$\tfrac3{10}$ & $1$\\
\noalign{\hrule}
\end{tabular}
\end{center}

\texttt{CUSPSANDWIDMAKE(N)} ---
   returns a set of inequivalent cusps for $\Gamma_1(N)$,
   and corresponding widths.
   Output has the form \texttt{[CUSPLIST,WIDTHLIST]}.

\exampleshade{
> CUSPSANDWIDMAKE1(10);
$$
\left[\left[
\mbox{oo}, 0, \frac12, \frac13, \frac14, \frac15, \frac25, \frac3{10}\right],
 [1, 10, 5, 10, 5, 2, 2, 1]\right]
$$
}

\subsection{Orders at cusps}
\mylabel{subsec:cuspords}

We will use Biagioli's \cite{Bi89} results for theta-functions to calculate orders at
cusps of generalized eta-products.
We define the theta-function
\begin{equation}
\theta_{\delta;g}(\tau) = q^{ (\delta-2g)^2/(8\delta)}
\prod_{m=1}^\infty (1 - q^{m\delta-g})(1-q^{m\delta-(\delta-g)})(1-q^{m\delta}),
\mylabel{eq:thetadgdef}
\end{equation}
for $0 < g < \delta$. This corresponds to Biagioli's function $f_{\delta,g}$ 
\cite[p.277]{Bi89}.  The classical Dedekind eta-function
can be written as
\beq
\eta(\tau) = \theta_{3;1}(\tau),
\mylabel{eq:etatheta}
\eeq
and
the generalized Dedekind eta-function can be written as
\begin{equation} 
\eta_{\delta;g}(\tau) = \frac{\theta_{\delta;g}(\tau)}{\eta(\delta\tau)}
 = \frac{\theta_{\delta;g}(\tau)}{\theta_{3\delta;\delta}(\tau)}.
\mylabel{eq:getatheta}
\end{equation}
Biagioli \cite{Bi89} has calculated the invariant order of $\theta_{\delta;g}(\tau)$
at any cusp. Using \eqn{getatheta} this gives a method for calculating the
invariant order at any cusp of a generalized eta-product.

\begin{theorem}[\cite{Bi89}(Lemma 3.2, p.285)]
\mylabel{thm:biaord}
The order at the cusp $s=\frac{b}{c}$ (assuming $(b,c)=1$)
of the theta function $\theta_{g;\delta}(\tau)$ (defined above
and assuming $\delta\nmid g$) is
\begin{equation}
\ord(\theta_{g;\delta}(\tau),s)
=\frac{e^2}{2\delta}\left(
                \frac{bg}{e} - \left[\frac{bg}{e}\right] - \frac{1}{2} \right)^2,
\mylabel{eq:Bord}
\end{equation}
where $e=(\delta,c)$ and $\left[\hphantom{\frac{a}{b}}\right]$
is the greatest integer function.
\end{theorem}

\texttt{Bord}$(\delta,g,a,c)$ --- returns the order of $\theta_{\delta;g}(\tau)$ at
the cusp $a/c$, assuming $(a,c)=1$ and $\delta\nmid g$.

\texttt{getacuspord}$(\delta,g,a,c)$ --- 
returns the order of the generalized eta-function $\eta_{\delta,g}(\tau)$ at
the cusp $a/c$, assuming $(a,c)=1$ and $\delta\nmid g$.

\exampleshade{
>  getacuspord(50,1,4,29);
$$
\frac{1}{600}
$$
}
We see that
$$
\ord \left(\eta_{50;1}(\tau), \frac{4}{29}\right) = \frac{1}{600}.
$$

   Let $G$ be a generalized eta-product corresponding to the    
   getalist $L$. The following MAPLE procedure calculates the invariant order 
   $\ord(G,\zeta)$ for any cusp $\zeta$.

\texttt{getaprodcuspord(L,cusp)} --- returns of the generalized eta-product
corresponding to the geta-list $L$ at the given cusp. The cusp is either
a rational or \texttt{oo} (infinity).

\exampleshade{
> GL:=[[4,1,16],[4,2,-4]];
$$
                    [[4, 1, 16], [4, 2, -4]]
$$
> getaprodcuspord(GL,1/2);
$$
                               -1
$$
}

We see that
$$
\ord\left({\frac{\eta_{4;1}(\tau)^{16}}{\eta_{4;2}(\tau)^{4}}},
\frac{1}{2}\right)
= -1.
$$

Following \cite[p.275]{Bi89}, \cite[p.91]{Ra} we consider the order of
a function $f$ with respect to a congruence subgroup $\Gamma$ at the
cusp $\zeta\in \mathbb{Q} \cup \{\infty\}$ and denote this by
\beq
\ORD(f,\zeta,\Gamma) = \kappa(\zeta,\Gamma)\,\ord(f,\zeta)
\mylabel{eq:ORD}.
\eeq


\texttt{getaprodcuspORDS}$(L,S,W)$        ---
returns a list of orders $\ORD(G,\zeta,\Gamma_1(N))$ where
$G$ 
is the  generalized eta-product corresponding to the    
   getalist $L$, $\zeta \in S$ (list of inequivalent cusps of 
$\Gamma_1(N)$) and $W$ is a list of corresponding fan-widths.

\exampleshade{
> CW4:=CUSPSANDWIDMAKE1(4);
$$
\left[\left[\infty,0,\frac{1}{2}\right],[1, 4, 1]\right]
$$
> GL:=[[4,1,16],[4,2,-4]];
$$
                    [[4, 1, 16], [4, 2, -4]]
$$
> getaprodcuspORDS(GL,CW4[1],CW4[2]);
$$
                           [0, 1, -1]
$$
}

We know that the generalized eta-product
$$
f(\tau) = \frac{\eta_{4;1}(\tau)^{16}}{\eta_{4;2}(\tau)^{4}}
$$
is a modular function on $\Gamma_1(4)$. We calculated
$\ORD(f,\zeta,\Gamma_1(4))$ at each cusp $\zeta$ of $\Gamma_1(4)$.

\begin{center}
\begin{tabular}{|c|c|}
\noalign{\hrule}
$\zeta$ & $\ORD(f,\zeta,\Gamma_1(4)$\\
\noalign{\hrule}
$\infty$ & $0$\\
\noalign{\hrule}
$0$ & $1$ \\
\noalign{\hrule}
$\frac{1}{2}$  & $-1$\\
\noalign{\hrule}
\end{tabular}
\end{center}
Observe that the total order of $f$ with respect to $\Gamma_1(4)$ is $0$:
$$
\ORD(f,\Gamma_1(4)) = \sum_{\zeta \in \mathcal{S}} \ORD(f,\zeta,\Gamma_1(4))= 0 + 1 - 1 = 0,
$$
in agreement with the valence formula. See Theorem \thm{val}
below. Here $\mathcal{S}$ is the set of inequivalent cusps of $\Gamma_1(4)$.

\subsection{Proving theta-function identities}
\mylabel{subsec:provethetaid}

Our method for proving theta-function or generalized eta-product identities depends
on
\begin{theorem}[The Valence Formula \cite{Ra}(p.98)]
\mylabel{thm:val}
Let $f\ne0$ be a modular form of weight $k$ with respect to a subgroup $\Gamma$ of finite index
in $\Gamma(1)=\SL_2(\mathbb{Z})$. Then
\beq
\ORD(f,\Gamma) = \frac{1}{12} \mu \, k,
\mylabel{eq:valform}
\eeq
where $\mu$ is the index of $\widehat{\Gamma}$ in $\widehat{\Gamma(1)}$,
$$
\ORD(f,\Gamma) := \sum_{\zeta\in R^{*}} \ORD(f,\zeta,\Gamma),
$$
$R^{*}$ is a fundamental region for $\Gamma$,
and
$\ORD(f,\zeta,\Gamma)$ is given in equation \eqn{ORD}.
\end{theorem}
\begin{remark}
For $\zeta\in\mathfrak{h}$,
$\ORD(f,\zeta,\Gamma)$ is defined in terms of 
 the invariant order $\ord(f,\zeta)$, which  is interpreted
in the usual sense. See \cite[p.91]{Ra} for details of this and the 
notation used.
\end{remark}

Since any generalized eta-product has weight $k=0$ and has no zeros and no
poles on the upper-half plane we have
\begin{cor}
\mylabel{cor:valcor}
Let $f_1(\tau)$, $f_2(\tau)$, \dots, $f_n(\tau)$ be generalized eta-products that
are modular functions on $\Gamma_1(N)$. Let $\mathcal{S}_N$ be a set of inequivalent
cusps for $\Gamma_1(N)$. Define the constant
\beq
B = \sum_{\substack{s\in\mathcal{S}_N\\s\ne \infty}}
        \mbox{min}
        (\left\{\ORD(f_j,s,\Gamma_1(N))\,:\, 1 \le j \le n\right\} \cup \{0\}),
\mylabel{eq:Bdef}
\eeq
and consider
\beq
g(\tau) := \alpha_1 f_1(\tau) + \alpha_2 f_2(\tau) + \cdots + \alpha_n f_n(\tau) + 1, 
\mylabel{eq:gdef}
\eeq
where each $\alpha_j\in\mathbb{C}$. Then
$$
g(\tau) \equiv 0
$$
if and only if
\beq
\ORD(g(\tau), \infty, \Gamma_1(N)) > -B.
\mylabel{eq:ORDBineq}
\eeq
\end{cor}

    To prove an alleged theta-function identity, we first rewrite it in the form
\begin{equation} 
    \alpha_1 f_1(\tau) + \alpha_2 f_2(\tau) + \cdots + \alpha_n f_n(\tau) + 1 = 0,
\mylabel{eq:fid}
\end{equation}
where each $\alpha_i\in\C$ and each $f_i(\tau)$ is a generalized eta-product of 
level $N$. We use the following algorithm:

        \vskip 10pt\noindent
{\it\footnotesize STEP 1}. \quad  Use Theorem \thm{gepmodfunc} to check that
$f_j(\tau)$ is a generalized eta-product on $\Gamma_1(N)$ for each
$1 \le j \le n$.

        \vskip 10pt\noindent
{\it\footnotesize STEP 2}. \quad  Use Theorem \thm{chokoopark} to
find a set $\mathcal{S}_N$ of inequivalent cusps for $\Gamma_1(N)$ and the
fan width of each cusp.

        \vskip 10pt\noindent
{\it\footnotesize STEP 3}. \quad  Use Theorem \thm{biaord} to
calculate the invariant order of each generalized eta-product
$f_j(\tau)$ at each cusp of $\Gamma_1(N)$.

        \vskip 10pt\noindent
{\it\footnotesize STEP 4}. \quad  Calculate
        $$
        B =
        \sum_{\substack{s\in\mathcal{S}_N\\s\ne \infty}}
        \mbox{min}
        (\left\{\ORD(f_j,s,\Gamma_1(N))\,:\, 1 \le j \le n\right\} \cup \{0\}).
        $$

        \vskip 10pt\noindent
{\it\footnotesize STEP 5}. \quad  Show that
        $$
        \ORD(g(\tau),\infty,\Gamma_1(N)) > -B
        $$
        where
        $$
        g(\tau) = \alpha_1 f_1(\tau) + \alpha_2 f_2(\tau) +
        \cdots + \alpha_n f_n(\tau) + 1.
        $$
        Corollary \corol{valcor} then implies that $g(\tau)\equiv0$ and
        hence the theta-function identity  \eqn{fid}.

To calculate the constant $B$ in \eqn{Bdef} and STEP 4 we use

\texttt{mintotORDS(L,n)} --- returns the constant $B$ in equation \eqn{Bdef}
where $L$ the array of ORDS:
$$
L := [\ORD(f_1), \ORD(f_2), \dots, \ORD(f_n)],
$$
where
$$
\ORD(f) = [\ORD(f,\zeta_1,\Gamma_1(N)), 
\ORD(f,\zeta_2,\Gamma_1(N)), \dots,
\ORD(f,\zeta_m,\Gamma_1(N))] 
$$
and $\zeta_1$, $\zeta_2$, \dots, $\zeta_m$ are the 
inequivalent cusps of $\Gamma_1(N)$.
Each $\ORD(f)$ is computed using \texttt{getaprodcuspORDS}.

\noindent
EXAMPLE:
As an example we prove Ramanujan's well-known identity
\beq
\prod_{n=1}^\infty \frac{(1-q^n)}{(1-q^{25n})} = R(q^5) - q - \frac{q^2}{R(q^5)},
\mylabel{eq:ramid5}
\eeq
where
$$
R(q) = \prod_{n=1}^\infty \frac{(1 - q^{5n-2})(1 - q^{5n-3})}
                               {(1 - q^{5n-1})(1 - q^{5n-4})}.
$$
We rewrite this identity as
\beq
\frac{\eta(\tau)}{\eta(25\tau)} =
\frac{\eta_{25;10}(\tau)}{\eta_{25;5}(\tau)} 
- 1 -
\frac{\eta_{25;5}(\tau)}{\eta_{25;10}(\tau)}.
\mylabel{eq:altramid5}
\eeq
Let 
\beq
g(\tau) = f_1(\tau) - f_2(\tau) + f_3(\tau) + 1,
\mylabel{eq:gdefeg}
\eeq
where
$$
f_1(\tau) =
\frac{\eta(\tau)}{\eta(25\tau)} =
\prod_{j=1}^{12} \eta_{25,j}(\tau),\quad
f_2(\tau)=\frac{\eta_{25;10}(\tau)}{\eta_{25;5}(\tau)},\quad
f_3(\tau) = \frac{1}{f_2(\tau)} =
\frac{\eta_{25;5}(\tau)}{\eta_{25;10}(\tau)}.
$$

\vskip 10pt\noindent
{\it\footnotesize STEP 1}. \quad  We check that each function is a modular function
on $\Gamma_1(25)$.
\mapleshade{
> f1:=mul(GETA(25,j), j=1..12): \\
> f2:=GETA(25,10)/GETA(25,5): \\
> f3:=1/f2: \\
> GP1:=GETAP2getalist(f1): \\
> GP2:=GETAP2getalist(f2): \\
> GP3:=GETAP2getalist(f3): \\
> Gamma1ModFunc(GP1,25),Gamma1ModFunc(GP2,25),Gamma1ModFunc(GP3,25); \\
$$
                            1, 1, 1
$$
}

\vskip 10pt\noindent
{\it\footnotesize STEP 2}. \quad  We find a set of inequivalent cusps for 
$\Gamma_1(25)$ and their fan widths.
\mapleshade{
> CW25:=CUSPSANDWIDMAKE1(25):\\
> cusps25:=CW25[1];
$$
[\mbox{oo}, 0, \frac{1}{2}, \frac{1}{3}, \frac{1}{4}, \frac{1}{5}, \frac{1}{6}, \frac{1}{7}, \frac{1}{8}, \frac{1}{9}, \frac{1}{10}, \frac{1}{11}, \frac{1}{12}, \frac{2}{5}, \frac{2}{25}, \frac{3}{5}, \frac{3}{10}, \frac{3}{25}, \frac{4}{5}, \frac{4}{25}, \frac{6}{25}, \frac{7}{10}, \frac{7}{25}, \frac{8}{25}, \frac{9}{10}, \frac{9}{25}, \frac{11}{25}, \frac{12}{25}]
$$
> widths25:=CW25[2];
$$
[1, 25, 25, 25, 25, 5, 25, 25, 25, 25, 5, 25, 25, 5, 1, 5, 5, 1, 
  5, 1, 1, 5, 1, 1, 5, 1, 1, 1]
$$
}

\vskip 10pt\noindent
{\it\footnotesize STEP 3}. \quad  We compute $\ORD(f_j,\zeta,\Gamma_1(25))$
for each $j$ and each cusp $\zeta$ of $\Gamma_1(25)$.
\mapleshade{
> ORDS1:=getaprodcuspORDS(GP1,cusps25,widths25);\\
$$
[-1, 1, 1, 1, 1, 0, 1, 1, 1, 1, 0, 1, 1, 0, -1, 0, 0, -1, 0, -1, 
  -1, 0, -1, -1, 0, -1, -1, -1]
$$
> ORDS2:=getaprodcuspORDS(GP2,cusps25,widths25);\\
$$
[-1, 0, 0, 0, 0, 0, 0, 0, 0, 0, 0, 0, 0, 0, 1, 0, 0, 1, 0, -1, 
  -1, 0, 1, 1, 0, -1, -1, 1]
$$
> ORDS3:=getaprodcuspORDS(GP3,cusps25,widths25);\\
$$
[1, 0, 0, 0, 0, 0, 0, 0, 0, 0, 0, 0, 0, 0, -1, 0, 0, -1, 0, 1, 1, 
  0, -1, -1, 0, 1, 1, -1]
$$
}

\vskip 10pt\noindent
{\it\footnotesize STEP 4}. \quad  We calculate the constant $B$ in \eqn{Bdef}.
\mapleshade{
> mintotORDS([ORDS1,ORDS2,ORDS3],3);\\
$$
                               -9
$$
}

\vskip 10pt\noindent
{\it\footnotesize STEP 5}. \quad  To prove the identity \eqn{ramid5} we need to
verify that 
$$
\ORD(g(\tau),\infty,\Gamma_1(25)) > 9.
$$
\mapleshade{
> JACL:=map(getalist2jacprod,[GP1,GP2,GP3]): \\
> JACID:=JACL[1]-JACL[2]+JACL[3]+1:\\
> QJ:=jac2series(JACID,100): \\
> series(QJ,q,100);
$$
\mbox{O}\left(q^{99}\right)
$$
}
This completes the proof of the identity \eqn{ramid5}.
We only had to show that the coefficient of 
$q^j$ was zero in the $q$-expansion of 
$g(\tau)$ for $j <\le 10$.  We actually did it for 
$j \le 98$ as a check.

STEPS 1--5 may be automated using

\texttt{provemodfuncid(JACID,N)} --- returns the constant $B$ in equation \eqn{Bdef}
and prints details of the verification and proof of the identity corresponding
to \texttt{JACID}, which is a linear combination of symbolic \texttt{JAC}-functions,
and $N$ is the level. If \texttt{xprint=true} then more details of the
verification are printed. When this function is called there is a query asking
whether to verify the identity. Enter \texttt{yes} to carry out the verification.

\exampleshade{
> provemodfuncid(JACID,25);\\
"TERM ", 1, "of ", 4, " *****************"\\
"TERM ", 2, "of ", 4, " *****************"\\
"TERM ", 3, "of ", 4, " *****************"\\
"TERM ", 4, "of ", 4, " *****************"\\
"mintotord = ", -9\\
"TO PROVE the identity we need to show that v[oo](ID) > ", 9\\
*** There were NO errors. \\
*** o Each term was modular function on\\
      Gamma1(25). \\
*** o We also checked that the total order of\\
      each term was zero.\\
*** o We also checked that the power of q was correct in\\
      each term.\\
"*** WARNING: some terms were constants. ***"\\
"See array CONTERMS."\\
To prove the identity we will need to verify if up to \\
q\carrot(9).\\
Do you want to prove the identity? (yes/no)\\
You entered yes.\\
We verify the identity to O(q\carrot(59)).\\
RESULT: The identity holds to O(q\carrot(59)).\\
CONCLUSION: This proves the identity since we had only\\
            to show that v[oo](ID) > 9.\\
}
                               9

\texttt{provemodfuncidBATCH(JACID,N)} --- is a version of 
\texttt{provemodfuncid} that prints less detail and does not query.

\exampleshade{
> provemodfuncidBATCH(JACID,25);\\
*** There were NO errors.  Each term was modular function on\\
    Gamma1(25). Also -mintotord=9. To prove the identity\\
    we need to  check up to O(q\carrot(11)).\\
    To be on the safe side we check up to O(q\carrot(59)).\\
*** The identity is PROVED!\\
}

\texttt{printJACIDORDStable()} --- prints an ORDs table for the $f_j$
and lower bound for $g$ after \texttt{provemodfuncid} is run.
Formatted output from our example is given below in Table \ref{tab:ordfs}. 
By summing the last column
we see that $B=-9$, which confirms an earlier calculation using 
\texttt{mintotORDS}.

\begin{table}[!htbp]
$$
\begin{array}{c|c|c|c|c}
\noalign{\hrule}
 \zeta&\ORD \left(f_1,\zeta\right) &\ORD \left(f_2,\zeta\right) &
\ORD \left(f_3,\zeta\right) &\mbox{Lower bound for $\ORD \left(g,\zeta\right)$} 
\\ 
\noalign{\hrule}
\frac{1}{2}&1&0&0&0\\ \frac{1}{3}&1&0&0&0
\\ \frac{1}{4}&1&0&0&0\\ \frac{1}{5}&0&0&0&0
\\ \frac{1}{6}&1&0&0&0\\ \frac{1}{7}&1&0&0&0
\\ \frac{1}{8}&1&0&0&0\\ \frac{1}{9}&1&0&0&0
\\ \frac{1}{10}&0&0&0&0\\ \frac{1}{11}&1&0&0&0
\\ \frac{1}{12}&1&0&0&0\\ \frac{2}{5}&0&0&0&0
\\ {\frac {2}{25}}&-1&1&-1&-1\\ 
\frac{3}{5}&0&0&0&0\\ \frac{3}{10}&0&0&0&0\\ {\frac 
{3}{25}}&-1&1&-1&-1\\ \frac{4}{5}&0&0&0&0
\\ {\frac {4}{25}}&-1&-1&1&-1\\ {
\frac {6}{25}}&-1&-1&1&-1\\ {\frac {7}{10}}&0&0&0&0
\\ {\frac {7}{25}}&-1&1&-1&-1\\ {
\frac {8}{25}}&-1&1&-1&-1\\ {\frac {9}{10}}&0&0&0&0
\\ {\frac {9}{25}}&-1&-1&1&-1\\ {
\frac {11}{25}}&-1&-1&1&-1\\ {\frac {12}{25}}&-1&1&-
1&-1\\
\noalign{\hrule}
\end{array}
$$
\caption{Orders at the cusps of $\Gamma_1(25)$ of the functions 
$f_1$, $f_2$, $f_3$ and $g$  in 
\eqn{gdefeg} needed in the proof of Ramanujan's identity \eqn{altramid5}.
This table was produced by
\texttt{printJACIDORDStable()}.}                                              
\mylabel{tab:ordfs}
\end{table}

\section{Generalized Ramanujan-Robins identities}
\mylabel{sec:genramarobins}
\setcounter{equation}{0}

As an application of our \texttt{thetaids} package we show how to find and prove
generalized eta-product identities due to Ramanujan and Robins, and some natural
extensions. In Section \sect{intro} we defined 
the
functions $G(n,N,\chi)$ and 
$H(n,N,\chi)$, where $\chi$ is a non-principal
real Dirichlet character mod $N$ satisfying
$\chi(-1)=1$. 
Robins \cite{Ro-thesis} proved the following striking analogue
of Ramanujan's identity \eqn{ENTRY3.4} (or \eqn{altENTRY3.4}):
\beq
G(3) \, H(1) - G(1)\,H(3) =1,
\mylabel{eq:Robins13}
\eeq
where
$$
G(n) = \frac{1}{\eta_{13;1,3,4}(n\tau)},
\quad
H(n) = \frac{1}{\eta_{13;2,5,6}(n\tau)}.
$$
Equation \eqn{Robins13} is a restatement of
\eqn{RR13b}. In this case $N=13$ 
and $\chi = \leg{\cdot}{13}$ is the Legendre
symbol.

We will also consider
\beq
G^{*}(n,N,\chi) = G^{*}(n) := 
\prod_{\substack{\chi(g) = 1\\ 0 < g < \frac{N}{2}}}
\frac{1}{\eta^{*}_{N;g}(n\tau)}, \quad
H^{*}(n,N,\chi) = H^{*}(n) := 
\prod_{\substack{\chi(g) = -1\\ 0 < g < \frac{N}{2}}}
\frac{1}{\eta^{*}_{N;g}(n\tau)},
\mylabel{eq:GHstardefs}
\eeq
where
\beq
\eta^{*}_{\delta;g}(\tau) =q^{\frac{\delta}{2} P_2(g/\delta) }
\prod_{m\equiv \pm g\pmod{\delta}} (1 - (-q)^m).
\mylabel{eq:starGeta}
\eeq
We note that
$$
\eta^{*}_{\delta;g} (\tau) = \omega_{\delta;g} \eta_{\delta;g}(\tau + \pi i),
$$
where $\omega_{\delta;g}$ is a root of unity.
Using the notation \eqn{GHdefs} (with $N=5$ and $\chi(\cdot) = \leg{\cdot}{5}$,
the Legendre symbol) we may rewrite Ramanujan's 
identities
\eqn{ENTRY3.1}, \eqn{ENTRY3.4} as
\begin{align*}
     G(1)^{11}H(1) - G(1)H(1)^{11} &= 1+11G(1)^6H(1)^6,\\
     H(1)G(11) - G(1)H(11) &= 1,                           
\end{align*}
respectively.

We have written a number of specialized functions for the purpose 
of finding and proving identities for these more general 
$G$- and $H$-functions.  We have collected these functions into 
the new \texttt{ramarobinsids} package.
Go to
\begin{center}
\url{http://qseries.org/fgarvan/qmaple/ramarobinsids}
\end{center}
and  follow the directions on that page. This package requires
both the \texttt{qseries} and \texttt{thetaids} packages.

\subsection{Some \maple functions}
\mylabel{subsec:maplefuncs}

\texttt{Geta(g,d,n)} ---  returns the generalized eta-function 
$\eta_{d;g}(n\tau)$ in symbolic \texttt{JAC}-form.

\texttt{GetaB(g,d,n)} ---  returns \texttt{Geta(g,d,n)} without the
the $q^{\frac{d}{2} P_2(g/d) }$ factor.

\texttt{GetaL(L,d,n)} ---  returns the generalized eta-product 
corresponding to the geta-list in \texttt{JAC}-form 
with $\tau$ replaced by $n\tau$.

\texttt{GetaBL(L,d,n)} ---  returns the generalized eta-product 
\texttt{GetaL(g,d,n)} without the $q$-factor.

\texttt{GetaEXP(g,d,n)} ---  returns lowest power of $q$ in  
$\eta_{d;g}(n\tau)$.                                

\texttt{GetaLEXP(L,d,n)} ---  returns lowest power of $q$ for the
generalized eta-product corresponding to 
\texttt{GetaL(L,d,n)}.

\texttt{MGeta(g,d,n)} --- $\eta^{*}$ analogue of \texttt{Geta(g,d,n)}

\texttt{MGetaL(L,d,n)} --- $\eta^{*}$ analogue of \texttt{GetaL(L,d,n)}

\texttt{Eeta(n)} --- returns Dedekind eta-function $\eta(n\tau)$ in 
\texttt{JAC}-form.

\exampleshade{
> with(ramarobinsids):\\
> Geta(1,5,2);
$$
\frac{ q^{1/30} JAC(2,10,\infty)}{JAC(0,10,\infty)}
$$
> GetaB(1,5,2);
$$
\frac{JAC(2,10,\infty)}{JAC(0,10,\infty)}
$$
> GetaEXP(1,5,2);
$$
\frac{1}{30}
$$
> GetaL([1,3,4],13,1);
$$
\frac{q^{1/4}}{JAC(0,13,\infty)^3}
JAC(1,13,\infty) JAC(3,13,\infty) JAC(4,13,\infty)
$$
> GetaLB([1,3,4],13,1);
$$
\frac{JAC(1,13,\infty) JAC(3,13,\infty) JAC(4,13,\infty)}
{JAC(0,13,\infty)^3}
$$
> GetaLEXP([1,3,4],13,1);
$$
\frac{1}{4}
$$
> MGeta(1,5,2);
$$
{\frac {{q}^{1/30}{\it JAC} \left( 2,10,\infty  \right) {\it JAC}
 \left( 4,40,\infty  \right)  \left( {\it JAC} \left( 0,20,\infty 
 \right)  \right) ^{2}}{{\it JAC} \left( 0,10,\infty  \right) {\it JAC
} \left( 0,40,\infty  \right)  \left( {\it JAC} \left( 2,20,\infty 
 \right)  \right) ^{2}}}
$$
> MGetaL([1,3,4],13,1);
$$
{\frac {\sqrt [4]{q}{\it JAC} \left( 1,13,\infty  \right) {\it JAC}
 \left( 2,52,\infty  \right) {\it JAC} \left( 0,26,\infty  \right) {
\it JAC} \left( 3,13,\infty  \right) {\it JAC} \left( 6,52,\infty 
 \right)  \left( {\it JAC} \left( 4,26,\infty  \right)  \right) ^{2}{
\it JAC} \left( 8,26,\infty  \right) }{{\it JAC} \left( 0,13,\infty 
 \right) {\it JAC} \left( 0,52,\infty  \right)  \left( {\it JAC}
 \left( 1,26,\infty  \right)  \right) ^{2} \left( {\it JAC} \left( 3,
26,\infty  \right)  \right) ^{2}{\it JAC} \left( 4,13,\infty  \right) 
{\it JAC} \left( 8,52,\infty  \right) }}
$$
> Eeta(3);
$$              
q^{1/8} JAC(0,3,\infty)
$$
}

\texttt{CHECKRAMIDF(SYMF,ACC,T)} --- checks whether a certain symbolic
expression of $G$- and $H$-functions is an eta-product. This assumes
that $G(n)$, $H(n)$, $GM(n)$, $HM(n)$ have already been defined.
$GM$ and $HM$ are the $\eta^{*}$ analogues of  $G$, $H$.
The \texttt{SYMF} symbolic form is written in terms of 
$\_G$, $\_H$, $\_GM$, $\_HM$. \texttt{ACC} is an upperbound on the absolute
value of exponents allowed in the formal product, \texttt{T} is highest
power of $q$ considered. This procedure returns a list of exponents in the formal
product if it is a likely eta-product otherwise it returns \texttt{NULL}.
A number of global variables are also assigned. The main ones are
\begin{itemize}
\item
\texttt{\_JFUNC}: \texttt{JAC}-expression of \texttt{SYMF}.
\item
\texttt{LQD}: lowest power of $q$.
\item
\texttt{RID}: the conjectured eta-product.
\item
\texttt{ebase}: base of the conjectured eta-product.
\item
\texttt{SYMID}: symbolic form of the identity
\end{itemize} 

\exampleshade{
> with(qseries):\\
> with(thetaids):\\
> with(ramarobinsids):\\
> G:=j->1/GetaL([1,3,4],13,j): H:=j->1/GetaL([2,5,6],13,j):\\
> GM:=j->1/MGetaL([1,3,4],13,j): HM:=j->1/MGetaL([2,5,6],13,j):\\
> GE:=j->-GetaLEXP([1,3,4],13,j): HE:=j->-GetaLEXP([2,5,6],13,j):\\
> GHID:=(\_G(1)*\_G(2)+\_H(1)*\_H(2))/(\_G(2)*\_H(1)-\_G(1)*\_H(2));
$$
GHID := {\frac {{\it \_G} \left( 1 \right) {\it \_G} \left( 2 \right) +{\it 
\_H} \left( 1 \right) {\it \_H} \left( 2 \right) }{{\it \_G} \left( 2
 \right) {\it \_H} \left( 1 \right) -{\it \_G} \left( 1 \right) {\it 
\_H} \left( 2 \right) }}
$$
> CHECKRAMIDF(GHID,10,50);
\begin{align*}
& [-2, 0, -2, 0, -2, 0, -2, 0, -2, 0, -2, 0, 0, 0, -2, 0, -2, 0,  \\
& \quad -2, 0, -2, 0, -2, 0, -2, 0, -2, 0, -2, 0, -2, 0, -2, 0, -2, 0, 
  -2, 0, 0]
\end{align*}
> ebase;
$$
                               26
$$
> \_JFUNC;
\begin{align*}
&(
-{q}^{3}{\it JAC} \left( 1,13,\infty  \right) 
{\it JAC} \left( 3,13,\infty  \right) 
{\it JAC} \left( 4,13,\infty  \right) 
{ \it JAC} \left( 2,26,\infty  \right) 
{\it JAC} \left( 6,26,\infty \right) 
{\it JAC} \left( 8,26,\infty  \right)  \\
& \quad -{\it JAC} \left( 2,13, \infty  \right) 
{\it JAC} \left( 5,13,\infty  \right) 
{\it JAC} \left( 6,13,\infty  \right) 
{\it JAC} \left( 4,26,\infty  \right) 
{ \it JAC} \left( 10,26,\infty  \right) 
{\it JAC} \left( 12,26,\infty \right)
) \\
& 
/
q \left( q{\it JAC} \left( 2,26,\infty  \right) 
{\it JAC} \left( 6,26,\infty  \right) 
{\it JAC} \left( 8,26,\infty  \right) 
{ \it JAC} \left( 2,13,\infty  \right) 
{\it JAC} \left( 5,13,\infty \right) 
{\it JAC} \left( 6,13,\infty  \right)  \right.\\
& \left. \quad -{\it JAC} \left( 1,13, \infty  \right) 
{\it JAC} \left( 3,13,\infty  \right) 
{\it JAC} \left( 4,13,\infty  \right) 
{\it JAC} \left( 4,26,\infty  \right) 
{ \it JAC} \left( 10,26,\infty  \right) 
{\it JAC} \left( 12,26,\infty \right)  \right)
\end{align*}
> LDQ;
$$
                               -1
$$
> RID;
$$
{\frac { \left( \eta \left( 13\,\tau \right)  \right) ^{2} \left( \eta
 \left( 2\,\tau \right)  \right) ^{2}}{ \left( \eta \left( 26\,\tau
 \right)  \right) ^{2} \left( \eta \left( \tau \right)  \right) ^{2}}}
$$
> SYMID;
$$
{\frac {{\it \_G} \left( 1 \right) {\it \_G} \left( 2 \right) +{\it 
\_H} \left( 1 \right) {\it \_H} \left( 2 \right) }{{\it \_G} \left( 2
 \right) {\it \_H} \left( 1 \right) -{\it \_G} \left( 1 \right) {\it 
\_H} \left( 2 \right) }}
={\frac { \left( \eta \left( 13\,\tau \right) 
 \right) ^{2} \left( \eta \left( 2\,\tau \right)  \right) ^{2}}{
 \left( \eta \left( 26\,\tau \right)  \right) ^{2} \left( \eta \left( 
\tau \right)  \right) ^{2}}}
$$
> etamake(jac2series(\_JFUNC,1001),q,1001);
$$
{\frac {  \eta \left( 13\,\tau \right)^{2} \eta \left( 2\,\tau \right)^{2}}
{\eta \left( 26\,\tau \right)^{2}\eta \left( \tau \right)^{2}}}
$$
}
It seems that 
\beq
{\frac {{G} \left( 1 \right) {G} \left( 2 \right) +
{H} \left( 1 \right) {H} \left( 2 \right) }{{G} \left( 2
 \right) {H} \left( 1 \right) -{G} \left( 1 \right) 
{H} \left( 2 \right) }}
=
{\frac {  \eta \left( 13\,\tau \right)^{2} \eta \left( 2\,\tau \right)^{2}}
{\eta \left( 26\,\tau \right)^{2}\eta \left( \tau \right)^{2}}}
\mylabel{eq:newRRID13}
\eeq
when $N=13$ and $\chi(\cdot) = \leg{\cdot}{13}$, at least up to $q^{1000}$.

\exampleshade{
> RRID1:=\_JFUNC-Eeta(13)\carrot2*Eeta(2)\carrot2/Eeta(26)\carrot2/Eeta(1)\carrot2:\\
> JRID1:=processjacid(RRID1):\\
> jmxperiod;
$$
                               26
$$
> provemodfuncidBATCH(JRID1,26);\\
*** There were NO errors.  Each term was modular function on\\
    Gamma1(26). Also -mintotord=18. To prove the identity\\
    we need to  check up to O(q\carrot(20)).\\
    To be on the safe side we check up to O(q\carrot(70)).\\
*** The identity is PROVED!
}
Thus identity \eqn{newRRID13} is proved.

The search for and proof of such identities may
 be automated.

\subsection{Ten types of identities for Ramanujan's functions $G(q)$ and $H(q)$}
\mylabel{subsec:10RGH}

We consider ten types of identities. We write a \maple function to 
search for and prove identities of each type. Here we assume $N=5$
and $\chi(\cdot) = \leg{\cdot}{5}$. We continue to use the notation \eqn{GHdefs}.

In this section
$$
G(1)=G(1,5,\leg{\cdot}{5}) = \frac{1}{\eta_{5;1}(\tau)} 
= \frac{q^{-1/60}}{(q,q^4;q^5)_\infty},
\quad
H(1)=H(1,5,\leg{\cdot}{5}) = \frac{1}{\eta_{5;2}(\tau)} 
= \frac{q^{11/60}}{(q^2,q^3;q^5)_\infty}.
$$

\exampleshade{
> with(qseries):\\
> with(thetaids):\\
> with(ramarobinsids):\\
> G:=j->1/GetaL([1],5,j): H:=j->1/GetaL([2],5,j):\\
> GM:=j->1/MGetaL([1],5,j): HM:=j->1/MGetaL([2],5,j):\\
> GE:=j->-GetaLEXP([1],5,j): HE:=j->-GetaLEXP([2],5,j):
}

\subsubsection{Type 1}
\mylabel{subsubsec:type5.1}

We consider identities of the form
$$
G(a)\,H(b) \pm G(b)\,H(a) = f(\tau),
$$
where $f(\tau)$ is an eta-product and $a$, $b$ are positive relatively prime
integers.

\texttt{findtype1(T)} --- cycles through symbolic expressions
$$
\_G(a)\,\_H(b) + c \,\_G(b)\,\_H(a)
$$
where $ 2\le n \le T$, $ab=n$, $(a,b)=1$,  $b < a$, $c\in\{-1,1\}$,  and
\beq
\texttt{GE}(a) + \texttt{HE}(b) - (\texttt{GE}(b) + \texttt{HE}(a)) 
=\frac{1}{5}\left(b - a\right)
\in
\mathbb{Z},
\mylabel{eq:ft1cond}
\eeq
using \texttt{CHECKRAMIDF}
to check whether the expression corresponds to a likely eta-product, and if so
uses \texttt{provemodfuncidBATCH} to prove it. Condition \eqn{ft1cond}
eliminates the case of fractional powers of $q$, which in our case means
$a \equiv b$ $\mod 5$. The procedure also returns a list of \texttt{[a,b,c]}
which give identities.

\exampleshade{
> proveit:=true:\\
> findtype1(11);\\
*** There were NO errors.  Each term was modular function on\\
    Gamma1(30). Also -mintotord=8. To prove the identity\\
    we need to  check up to O(q\carrot(10)).\\
    To be on the safe side we check up to O(q\carrot(68)).\\
*** The identity below is PROVED!\\
\quad [6, 1, -1]
$$
       \_G(6) \_H(1) - \_G(1) \_H(6) = 
\frac{\eta(6\tau) \eta(\tau)}{\eta(3\tau) \eta(2\tau)}
$$
*** There were NO errors.  Each term was modular function on\\
    Gamma1(55). Also -mintotord=40. To prove the identity\\
    we need to  check up to O(q\carrot(42)).\\
    To be on the safe side we check up to O(q\carrot(150)).\\
*** The identity below is PROVED!\\
\quad [11, 1, -1]
\begin{gather*}
                \_G(11) \_H(1) - \_G(1) \_H(11) = 1\\
                   [[6, 1, -1], [11, 1, -1]]
\end{gather*}
}

\mapleshade{
> myramtype1 :=findtype1(36); 
$$
myramtype1 := [[6, 1, -1], [11, 1, -1], [7, 2, -1], [16, 1, -1], [8, 3, -1], [9, 4, -1], [36, 1, -1]]
$$
}
This also produced the following identities with proofs (some output
omitted):

\begin{alignat}{3}
G(6)\,H(1) - G(1)\,H(6) &= 
\frac{\eta(\tau)\eta(6\tau)}{\eta(2\tau)\eta(3\tau)},&\quad &\Gamma_1(30), &\quad -B=8,
\mylabel{eq:RR511}\\
G(11)\,H(1) - G(1)\,H(11) &= 
1
,&\quad &\Gamma_1(55), &\quad -B=40,
\mylabel{eq:RR512}\\
G(7)\,H(2) - G(2)\,H(7) &= 
\frac{\eta(\tau)\eta(14\tau)}{\eta(2\tau)\eta(7\tau)},&\quad &\Gamma_1(70), &\quad -B=48,
\mylabel{eq:RR513}\\
G(16)\,H(1) - G(1)\,H(16) &= 
\frac{\eta(4\tau)^{2}}{\eta(2\tau)\eta(8\tau)},&\quad &\Gamma_1(80), &\quad -B=64,
\mylabel{eq:RR514}\\
G(8)\,H(3) - G(3)\,H(8) &= 
\frac{\eta(\tau)\eta(4\tau)\eta(6\tau)\eta(24\tau)}{\eta(2\tau)\eta(3\tau)\eta(8\tau)\eta(12\tau)},&\quad &\Gamma_1(120), &\quad -B=128,
\mylabel{eq:RR515}\\
G(9)\,H(4) - G(4)\,H(9) &= 
\frac{\eta(\tau)\eta(6\tau)^{2}\eta(36\tau)}{\eta(2\tau)\eta(3\tau)\eta(12\tau)\eta(18\tau)},&\quad &\Gamma_1(180), &\quad -B=288,
\mylabel{eq:RR516}\\
G(36)\,H(1) - G(1)\,H(36) &= 
\frac{\eta(4\tau)\eta(6\tau)^{2}\eta(9\tau)}{\eta(2\tau)\eta(3\tau)\eta(12\tau)\eta(18\tau)},&\quad &\Gamma_1(180), &\quad -B=288.
\mylabel{eq:RR517}
\end{alignat}

We have included the relevant groups $\Gamma_1(N)$ and values of $B$ (see \eqn{Bdef}
and \eqn{ORDBineq}).
These identities are known and are equations (3.9), (3.5), (3.10), (3.6), (3.12), 
(3.14), and (3.15) in \cite{Be-etal07} respectively.

\subsubsection{Type 2}
\mylabel{subsubsec:type5.2}

We consider identities of the form
$$
G(a)\,G(b) \pm H(a)\,H(b) = f(\tau),
$$
where $f(\tau)$ is an eta-product and $a$, $b$ are positive relatively prime
integers.

\texttt{findtype2(T)} --- cycles through symbolic expressions
$$
\_G(a)\,\_G(b) + c \,\_H(a)\,\_H(b)
$$
where $ 2\le n \le T$, $ab=n$, $(a,b)=1$,  $a \le b$, $c\in\{-1,1\}$,  and
\beq
\texttt{GE}(a) + \texttt{GE}(b) - (\texttt{HE}(a) + \texttt{HE}(b)) 
=-\frac{1}{5}\left(a + b\right)
\in
\mathbb{Z},
\mylabel{eq:ft2cond}
\eeq
using \texttt{CHECKRAMIDF}
to check whether the expression corresponds to a likely eta-product, and if so
uses \texttt{provemodfuncidBATCH} to prove it. Condition \eqn{ft2cond}
eliminates the case of fractional powers of $q$, which our case means
$a \equiv -b$ $\mod 5$. The procedure also returns a list of \texttt{[a,b,c]}
which give identities.

\mapleshade{
> findtype2(24); 
$$
[[1, 4, -1], [1, 4, 1], [2, 3, 1], [1, 9, 1], [1, 14, 1], [1, 24, 1]]
$$
}
This also produces the following identities with proofs:

\begin{alignat}{3}
G(1)\,G(4) - H(1)\,H(4) &= 
\frac{\eta(10\tau)^{5}}{\eta(2\tau)\eta(5\tau)^{2}\eta(20\tau)^{2}},&\quad &\Gamma_1(20), &\quad -B=4,\mylabel{eq:RR522}\\
G(1)\,G(4) + H(1)\,H(4) &= 
\frac{\eta(2\tau)^{4}}{\eta(\tau)^{2}\eta(4\tau)^{2}},&\quad &\Gamma_1(20), &\quad -B=4,\mylabel{eq:RR523}\\
G(2)\,G(3) + H(2)\,H(3) &= 
\frac{\eta(2\tau)\eta(3\tau)}{\eta(\tau)\eta(6\tau)},&\quad &\Gamma_1(30), &\quad -B=8,\mylabel{eq:RR524}\\
G(1)\,G(9) + H(1)\,H(9) &= 
\frac{\eta(3\tau)^{2}}{\eta(\tau)\eta(9\tau)},&\quad &\Gamma_1(45), &\quad -B=24,\mylabel{eq:RR525}\\
G(1)\,G(14) + H(1)\,H(14) &= 
\frac{\eta(2\tau)\eta(7\tau)}{\eta(\tau)\eta(14\tau)},&\quad &\Gamma_1(70), &\quad -B=48,\mylabel{eq:RR526}\\
G(1)\,G(24) + H(1)\,H(24) &= 
\frac{\eta(2\tau)\eta(3\tau)\eta(8\tau)\eta(12\tau)}{\eta(\tau)\eta(4\tau)\eta(6\tau)\eta(24\tau)},&\quad &\Gamma_1(120), &\quad -B=128,\mylabel{eq:RR527}
\end{alignat}
These identities are known and are equations 
(3.4), (3.3), (3.8), (3.7), (3.11), and (3.13) in \cite{Be-etal07} respectively.

\subsubsection{Type 3}
\mylabel{subsubsec:type5.3}

We consider identities of the form
$$
\frac{G(a_1)\,G(b_1) \pm H(a_1)\,H(b_1)}
{G(a_2)\,H(b_2) \pm H(a_2)\,G(b_2)}
 = f(\tau),
$$
which are not a quotient of Type 1 and 2 identities, and
where $f(\tau)$ is an eta-product, $a_1$, $b_1$ ,$a_2$, $b_2$ 
are positive relatively prime integers, and
$a_1 b_1 = a_2 b_2$.

\texttt{findtype3(T)} --- cycles through symbolic expressions
$$
\frac{\_G(a_1)\,\_G(b_1) + c_1 \,\_H(a_1)\,\_H(b_1)}
{\_G(a_2)\,\_H(b_2) +c_2\, \_H(a_2)\,\_G(b_2)}
$$
where $ 2\le n \le T$, $a_1b_1=a_2b_2=n$, $(a_1,b_1,a_2,b_2)=1$,  $a_1 \le b_1$, 
$b_2 < a_2$,
$c_1,c_2\in\{-1,1\}$,  and
\beq
\texttt{GE}(a_1) + \texttt{GE}(b_1) - (\texttt{HE}(a_1) + \texttt{HE}(b_1)),\quad
\texttt{GE}(a_2) + \texttt{HE}(b_2) - (\texttt{HE}(a_2) + \texttt{GE}(b_2)),\quad
\in
\mathbb{Z},
\mylabel{eq:ft3cond}
\eeq
and \texttt{[$a_2,b_2,c_2$]} is not an element of the list \texttt{myramtype1}
(found earlier by \texttt{findtype1}),
using \texttt{CHECKRAMIDF}
to check whether the expression corresponds to a likely eta-product, and if so
uses \texttt{provemodfuncidBATCH} to prove it. 
The procedure also returns lists \texttt{[$a_1,b_1,c_1,a_2,b_2,c_2$]}
which correspond to identities.

\mapleshade{
> findtype3(126); 
\begin{gather*}
[[3, 7, 1, 21, 1, -1], [2, 13, 1, 26, 1, -1], [1, 34, 1, 17, 2, -1], [1, 39, 1, 13, 3, -1], [1, 54, 1, 27, 2, -1], \\ 
[7, 8, 1, 56, 1, -1],
[3, 22, 1, 11, 6, -1], [2, 33, 1, 66, 1, -1], [4, 21, 1, 12, 7, -1], [1, 84, 1, 28, 3, -1], \\
[3, 32, 1, 96, 1, -1], [7, 18, 1, 14, 9, -1], [2, 63, 1, 126, 1, -1]]
\end{gather*}
}
This also produces the following identities with proofs:

\begin{alignat}{3}
\frac{G(3)\,G(7) + H(3)\,H(7)}{G(21)\,H(1) - H(21)\,G(1)} &= 
1
,&\quad &\Gamma_1(105), &\quad -B=192,
\mylabel{eq:RR531}\\
\frac{G(2)\,G(13) + H(2)\,H(13)}{G(26)\,H(1) - H(26)\,G(1)} &= 
1
,&\quad &\Gamma_1(130), &\quad -B=240,
\mylabel{eq:RR532}\\
\frac{G(1)\,G(34) + H(1)\,H(34)}{G(17)\,H(2) - H(17)\,G(2)} &= 
\frac{\eta(2\tau)\eta(17\tau)}{\eta(\tau)\eta(34\tau)},&\quad &\Gamma_1(170), &\quad -B=448,
\mylabel{eq:RR533}\\
\frac{G(1)\,G(39) + H(1)\,H(39)}{G(13)\,H(3) - H(13)\,G(3)} &= 
\frac{\eta(3\tau)\eta(13\tau)}{\eta(\tau)\eta(39\tau)},&\quad &\Gamma_1(195), &\quad -B=768,
\mylabel{eq:RR534}\\
\frac{G(1)\,G(54) + H(1)\,H(54)}{G(27)\,H(2) - H(27)\,G(2)} &= 
\frac{\eta(2\tau)\eta(3\tau)\eta(18\tau)\eta(27\tau)}{\eta(\tau)\eta(6\tau)\eta(9\tau)\eta(54\tau)},&\quad &\Gamma_1(270), &\quad -B=1008,
\mylabel{eq:RR535}\\
\star\quad\frac{G(7)\,G(8) + H(7)\,H(8)}{G(56)\,H(1) - H(56)\,G(1)} &= 
\frac{\eta(2\tau)\eta(28\tau)}{\eta(4\tau)\eta(14\tau)},&\quad &\Gamma_1(280), &\quad -B=1152,
\mylabel{eq:RR536}\\
\frac{G(3)\,G(22) + H(3)\,H(22)}{G(11)\,H(6) - H(11)\,G(6)} &= 
\frac{\eta(2\tau)\eta(33\tau)}{\eta(\tau)\eta(66\tau)},&\quad &\Gamma_1(330), &\quad -B=1600,
\mylabel{eq:RR537}\\
\frac{G(2)\,G(33) + H(2)\,H(33)}{G(66)\,H(1) - H(66)\,G(1)} &= 
\frac{\eta(3\tau)\eta(22\tau)}{\eta(6\tau)\eta(11\tau)},&\quad &\Gamma_1(330), &\quad -B=1600,
\mylabel{eq:RR538}\\
\star\quad\frac{G(4)\,G(21) + H(4)\,H(21)}{G(12)\,H(7) - H(12)\,G(7)} &= 
\frac{\eta(2\tau)\eta(3\tau)\eta(7\tau)\eta(12\tau)\eta(28\tau)\eta(42\tau)}{\eta(\tau)\eta(4\tau)\eta(6\tau)\eta(14\tau)\eta(21\tau)\eta(84\tau)},&\quad &\Gamma_1(420), &\quad -B=2688,
\mylabel{eq:RR539}\\
\star\quad\frac{G(1)\,G(84) + H(1)\,H(84)}{G(28)\,H(3) - H(28)\,G(3)} &= 
\frac{\eta(2\tau)\eta(3\tau)\eta(7\tau)\eta(12\tau)\eta(28\tau)\eta(42\tau)}{\eta(\tau)\eta(4\tau)\eta(6\tau)\eta(14\tau)\eta(21\tau)\eta(84\tau)},&\quad &\Gamma_1(420), &\quad -B=2688,
\mylabel{eq:RR5310}\\
\star\quad\frac{G(3)\,G(32) + H(3)\,H(32)}{G(96)\,H(1) - H(96)\,G(1)} &= 
\frac{\eta(2\tau)\eta(8\tau)\eta(12\tau)\eta(48\tau)}{\eta(4\tau)\eta(6\tau)\eta(16\tau)\eta(24\tau)},&\quad &\Gamma_1(480), &\quad -B=3072,
\mylabel{eq:RR5311}\\
\star\quad\frac{G(7)\,G(18) + H(7)\,H(18)}{G(14)\,H(9) - H(14)\,G(9)} &= 
\frac{\eta(2\tau)\eta(3\tau)\eta(42\tau)\eta(63\tau)}{\eta(\tau)\eta(6\tau)\eta(21\tau)\eta(126\tau)},&\quad &\Gamma_1(630), &\quad -B=5760,
\mylabel{eq:RR5312}\\
\frac{G(2)\,G(63) + H(2)\,H(63)}{G(126)\,H(1) - H(126)\,G(1)} &= 
\frac{\eta(3\tau)\eta(7\tau)\eta(18\tau)\eta(42\tau)}{\eta(6\tau)\eta(9\tau)\eta(14\tau)\eta(21\tau)},&\quad &\Gamma_1(630), &\quad -B=5760.
\mylabel{eq:RR5313}
\end{alignat}
The equations marked $\star$ appear to be new. The other equations
correspond to 
(3.16), (3.18), (3.35), (3.22), (3.41), (3.40) and (3.39) in \cite{Be-etal07},
and (1.24) in \cite{Ro-thesis}  respectively. 
We have corrected the statement of equation \cite[(1.24)]{Ro-thesis}.

\subsubsection{Type 4}
\mylabel{subsubsec:type5.4}

We consider identities of the form
$$
G^{*}(a)\,H^{*}(b) \pm G^{*}(b)\,H^{*}(a) = f(\tau),
$$
where $f(\tau)$ is an eta-product, $a$, $b$ are positive relatively prime
integers, and at least one of $a$, $b$ is even.

\texttt{findtype4(T)} --- cycles through symbolic expressions
$$
\_GM(a)\,\_HM(b) + c \,\_GM(b)\,\_HM(a)
$$
where $ 2\le n \le T$, $ab=n$, $(a,b)=1$,  $b < a$, $c\in\{-1,1\}$, 
\beq
\texttt{GE}(a) + \texttt{HE}(b) - (\texttt{GE}(b) + \texttt{HE}(a)) 
\in
\mathbb{Z},
\mylabel{eq:ft4cond}
\eeq
and least one of $a$, $b$ is even,
using \texttt{CHECKRAMIDF}
to check whether the expression corresponds to a likely eta-product, and if so
uses \texttt{provemodfuncidBATCH} to prove it. 
The procedure also returns a list of \texttt{[a,b,c]}
which give identities.

\mapleshade{
> findtype4(24); 
$$
[[6, 1, -1]]
$$
}
This also produces the following identity with proof:

\begin{equation}
G^{*}(6)\,H^{*}(1) - G^{*}(1)\,H^{*}(6) = 
\frac{\eta(\tau)\eta(4\tau)^{3}\eta(6\tau)^{3}\eta(24\tau)}{\eta(2\tau)^{3}\eta(3\tau)\eta(8\tau)\eta(12\tau)^{3}},\quad \Gamma_1(120), \quad -B=128.
\mylabel{eq:RR541}
\end{equation}
This corresponds to equation (3.28) in \cite{Be-etal07}.

\subsubsection{Type 5}
\mylabel{subsubsec:type5.5}

We consider identities of the form
$$
G^{*}(a)\,G^{*}(b) \pm H^{*}(a)\,H^{*}(b) = f(\tau),
$$
where $f(\tau)$ is an eta-product, $a$, $b$ are positive relatively prime
integers, and at least one of $a$, $b$ is even.

\texttt{findtype5(T)} --- cycles through symbolic expressions
$$
\_GM(a)\,\_GM(b) + c \,\_HM(a)\,\_HM(b)
$$
where $ 2\le n \le T$, $ab=n$, $(a,b)=1$,  $a \le b$, $c\in\{-1,1\}$, 
\beq
\texttt{GE}(a) + \texttt{GE}(b) - (\texttt{HE}(a) + \texttt{HE}(b)) 
\in
\mathbb{Z},
\mylabel{eq:ft5cond}
\eeq
and least one of $a$, $b$ is even,
using \texttt{CHECKRAMIDF}
to check whether the expression corresponds to a likely eta-product, and if so
uses \texttt{provemodfuncidBATCH} to prove it. 
The procedure also returns a list of \texttt{[a,b,c]}
which give identities.

\mapleshade{
> findtype5(24); 
$$
[[1, 4, 1], [2, 3, 1]]
$$
}
This also produces the following identity with proof:

\begin{alignat}{3}
G^{*}(1)\,G^{*}(4) + H^{*}(1)\,H^{*}(4) = 
\frac{\eta(4\tau)^{2}}{\eta(2\tau)\eta(8\tau)},\quad \Gamma_1(80), \quad -B=64,
\mylabel{eq:RR551}\\
G^{*}(2)\,G^{*}(3) + H^{*}(2)\,H^{*}(3) = 
\frac{\eta(2\tau)^{3}\eta(3\tau)\eta(8\tau)\eta(12\tau)^{3}}{\eta(\tau)\eta(4\tau)^{3}\eta(6\tau)^{3}\eta(24\tau)},\quad \Gamma_1(120), \quad -B=128.
\mylabel{eq:RR552}
\end{alignat}
These correspond to equations (3.26) and (3.27) in \cite{Be-etal07}.

\subsubsection{Type 6}
\mylabel{subsubsec:type5.6}

We consider identities of the form
$$
G(a)\,H^{*}(b) \pm G^{*}(a)\,H(b) = f(\tau),
$$
where $f(\tau)$ is an eta-product, $a$, $b$ are positive relatively prime
integers.

\texttt{findtype6(T)} --- cycles through symbolic expressions
$$
\_G(a)\,\_HM(b) + c \,\_GM(a)\,\_H(b)
$$
where $ 2\le n \le T$, $ab=n$, $(a,b)=1$,  $a \ge b$, $c\in\{-1,1\}$, 
using \texttt{CHECKRAMIDF}
to check whether the expression corresponds to a likely eta-product, and if so
uses \texttt{provemodfuncidBATCH} to prove it. 
The procedure also returns a list of \texttt{[a,b,c]}
which give identities.

\mapleshade{
> findtype6(24);
$$
                    [[1, 1, -1], [1, 1, 1]]
$$
}

This also produces the following identities with proof:
\begin{alignat}{3}
G(1)\,H^{*}(1) - G^{*}(1)\,H(1) &= 
{2}\,\frac{\eta(20\tau)^{2}}{\eta(2\tau)\eta(10\tau)},&\quad &\Gamma_1(20), &\quad -B=4,
\mylabel{eq:RR561}\\
G(1)\,H^{*}(1) + G^{*}(1)\,H(1) &= 
{2}\,\frac{\eta(4\tau)^{2}}{\eta(2\tau)^{2}},&\quad &\Gamma_1(20), &\quad -B=4.
\mylabel{eq:RR562}
\end{alignat}
These are equivalent to equations (3.25) and (3.24) in \cite{Be-etal07}.

\subsubsection{Type 7}
\mylabel{subsubsec:type5.7}
We consider identities of the form
$$
G^{*}(a)\,G(b) \pm H^{*}(a)\,H(b) = f(\tau),
$$
where $f(\tau)$ is an eta-product, $a$, $b$ are positive relatively prime
integers.

\texttt{findtype7(T)} --- cycles through symbolic expressions
$$
\_GM(a)\,\_G(b) + c \,\_HM(a)\,\_H(b)
$$
where $ 2\le n \le T$, $ab=n$, $(a,b)=1$,  $a \le b$, $c\in\{-1,1\}$, 
and both $a$, $b$ are odd,
using \texttt{CHECKRAMIDF}
to check whether the expression corresponds to a likely eta-product, and if so
uses \texttt{provemodfuncidBATCH} to prove it. 
The procedure also returns a list of \texttt{[a,b,c]}
which give identities.

\mapleshade{
> findtype7(24);
$$
                          [[1, 9, -1]]
$$
}

This also produces the following identity with proof:
\begin{equation}
G^{*}(1)\,G(9) - H^{*}(1)\,H(9) = 
\frac{\eta(\tau)\eta(12\tau)\eta(18\tau)^{2}}{\eta(2\tau)\eta(6\tau)\eta(9\tau)\eta(36\tau)},\quad \Gamma_1(180), \quad -B=288.
\mylabel{eq:RR571}
\end{equation}
This corresponds to (3.29) in \cite{Be-etal07}.

\subsubsection{Type 8}
\mylabel{subsubsec:type5.8}
We consider identities of the form
$$
G(1)^a\,H(a) \pm H(1)^a\,G(a) = f(\tau),
$$
where $f(\tau)$ is an eta-product, and $a >1$ is an integer.

\texttt{findtype8(T)} --- cycles through symbolic expressions
$$
\_G(1)^a\,\_H(a) + c \,\_H(1)^a\,\_G(a)
$$
where $ 2\le a \le T$, and  $c\in\{-1,1\}$,
using \texttt{CHECKRAMIDF}
to check whether the expression corresponds to a likely eta-product, and if so
uses \texttt{provemodfuncidBATCH} to prove it. 
The procedure also returns a list of \texttt{[a,c]}
which give identities.

\mapleshade{
> findtype8(24);
$$
                           [[3, -1]]
$$
}
This also produces the following identity with proof:
\begin{equation}
G(1)^{3}\,H(3) - H(1)^{3}\,G(3) = 
{3}\,\frac{\eta(15\tau)^{3}}{\eta(\tau)\eta(3\tau)\eta(5\tau)},\quad \Gamma_1(15), \quad -B=4.
\mylabel{eq:RR581}
\end{equation}
This is equivalent to equation (1.27) in Robin's thesis \cite{Ro-thesis}.

\subsubsection{Type 9}
\mylabel{subsubsec:type5.9}
We consider identities of the form
$$
G(1)^a\,H(1)^b - H(1)^a\,G(1)^b  + x = f(\tau),
$$
where $f(\tau)$ is an eta-product, and $a$, $b$ are positive integers,
and $x=0$ or $x=-1$.

\texttt{findtype9()} --- determines whether                     
$$
\_G(1)^a\,\_H(1)^b - \_H(1)^a\,\_G(1)^b + x
$$
is a likely eta-product for $x=0$ or $x=-1$ with $a$, $b$ smallest such
positive integers, 
using \texttt{CHECKRAMIDF}, and if so
uses \texttt{provemodfuncidBATCH} to prove it. 
The procedure also returns a list of \texttt{[a,b,x]}
which give identities.

\mapleshade{
> findtype9();
$$
                          [[11, 1, 1]]
$$
}
This also produces the following identity with proof:
\begin{equation}
G(1)^{11}\,H(1)^{1} - H(1)^{11}\,G(1)^{1} - 1 = 
{11}\,\frac{\eta(5\tau)^{6}}{\eta(\tau)^{6}},\quad \Gamma_1(5), \quad -B=2.
\mylabel{eq:RR591}
\end{equation}
This is equation (3.1) in \cite{Be-etal07}.

\subsubsection{Type 10}
\mylabel{subsubsec:type5.10}
We consider identities of the form
$$
\frac{G(a_1)\,H(b_1) \pm \,H(a_1)\,G(b_1)}
{G(a_2)\,H^{*}(b_2) \pm \, H(a_2)\,G^{*}(b_2)}
 = f(\tau),
$$
in which the numerator is not a Type 1 identity,
and
where $f(\tau)$ is an eta-product, $a_1$, $b_1$ ,$a_2$, $b_2$ 
are positive relatively prime integers, and $a_1 b_1 = a_2 b_2$.

\texttt{findtype10(T)} --- cycles through symbolic expressions
$$
\frac{\_G(a_1)\,\_H(b_1) + c_1 \,\_H(a_1)\,\_G(b_1)}
{\_G(a_2)\,\_HM(b_2) +c_2\, \_H(a_2)\,\_GM(b_2)}
$$
where $ 2\le n \le T$, $a_1b_1=a_2b_2=n$, $(a_1,b_1,a_2,b_2)=1$,  $a_1 >   b_1$, 
$b_2 < a_2$,
$c_1,c_2\in\{-1,1\}$,  and
\beq
\texttt{GE}(a_1) + \texttt{HE}(b_1) - (\texttt{HE}(a_1) + \texttt{GE}(b_1)),\quad
\texttt{GE}(a_2) + \texttt{HE}(b_2) - (\texttt{HE}(a_2) + \texttt{GE}(b_2)),\quad
\in
\mathbb{Z},
\mylabel{eq:ft10cond}
\eeq
and \texttt{[$a_1,b_1,c_1$]} is not an element of  the list \texttt{myramtype1}
(found earlier by \texttt{findtype1}),
using \texttt{CHECKRAMIDF}
to check whether the expression corresponds to a likely eta-product, and if so
uses \texttt{provemodfuncidBATCH} to prove it. 
The procedure also returns a list of \texttt{[$a_1,b_1,c_1,a_2,b_2,c_2$]}
which give identities.

\mapleshade{
> qthreshold:=3000:\\
> findtype10(120);
$$
        [[19, 4, -1, 76, 1, 1], [28, 3, -1, 12, 7, 1], 
          [12, 7, -1, 28, 3, 1]]
$$
}
This also produces the following identities with proof:
\begin{alignat}{3}
\frac{G(19)\,H(4) - H(19)\,G(4)}{G(76)\,H^{*}(1) + H(76)\,G^{*}(1)} &= 
\frac{\eta(2\tau)\eta(76\tau)}{\eta(4\tau)\eta(38\tau)},&\quad &\Gamma_1(380), &\quad -B=2160,
\mylabel{eq:RR5101}\\
\star\quad\frac{G(28)\,H(3) - H(28)\,G(3)}{G(12)\,H^{*}(7) + H(12)\,G^{*}(7)} &= 
\frac{\eta(\tau)\eta(4\tau)\eta(6\tau)\eta(14\tau)^{2}\eta(21\tau)}{\eta(2\tau)^{2}\eta(3\tau)\eta(7\tau)\eta(28\tau)\eta(42\tau)},&\quad &\Gamma_1(420), &\quad -B=2400,
\mylabel{eq:RR5102}\\
\star\quad\frac{G(12)\,H(7) - H(12)\,G(7)}{G(28)\,H^{*}(3) + H(28)\,G^{*}(3)} &= 
\frac{\eta(\tau)\eta(6\tau)^{2}\eta(14\tau)\eta(21\tau)\eta(84\tau)}{\eta(2\tau)\eta(3\tau)\eta(7\tau)\eta(12\tau)\eta(42\tau)^{2}},&\quad &\Gamma_1(420), &\quad -B=2400.
\mylabel{eq:RR5103}
\end{alignat}
Equation \eqn{RR5101} is (3.38) in \cite{Be-etal07}. The other type 10 
identities
appear to be new.

\section{More Generalized Ramanujan-Robins identities}
\mylabel{sec:moregenramarobins}
\setcounter{equation}{0}
We consider generalized Ramanujan-Robins identities associated with non-principal
real Dirichlet
characters $\chi$ mod $N$ for $N\le 60$, that satisfy $\chi(-1)=1$.
We found David Ireland's \textit{Dirichlet Character Table Generator}
\cite{DI-dctg} useful. See the website
\begin{center}
\url{http://www.di-mgt.com.au/dirichlet-character-generator.html}
\end{center}

\subsection{Mod $8$}
\mylabel{subsec:mod8}
There is only one non-principal character mod $8$ that satisfies $\chi(-1)=1$,
namely $\chi(\cdot)= \leg{8}{\cdot}$.
Here $\leg{8}{\cdot}$ is the Kronecker symbol.
In this section
$$
G(1)=G(1,8,\leg{8}{\cdot}) = \frac{1}{\eta_{8;1}(\tau)} 
= \frac{q^{-11/48}}{(q,q^7;q^8)_\infty},
\quad
H(1)=H(1,8,\leg{8}{\cdot}) = \frac{1}{\eta_{8;3}(\tau)} 
= \frac{q^{13/48}}{(q^3,q^5;q^8)_\infty}.
$$
These functions were considered by Robins \cite[pp.16-17]{Ro-thesis}.
They are also related to the G\"ollnitz-Gordon functions \cite{Go1965}, 
\cite{Go1967}:
\begin{align*}
S(q) &= \sum_{n=0}^\infty \frac{(-q;q^2)_n}{(q^2;q^2)_n}\,q^{n^2}
= \frac{1}{(q,q^4,q^7;q^8)_\infty},\\
T(q) &= \sum_{n=0}^\infty \frac{(-q;q^2)_n}{(q^2;q^2)_n}\,q^{n^2+2n}
= \frac{1}{(q^3,q^4,q^5;q^8)_\infty}.
\end{align*}
The ratio of these two functions is the famous Ramanujan-G\"ollnitz-Gordon continued
fraction \cite[Eq.(9.3)]{Be-overview01}
\begin{align*}
\frac{S(q)}{T(q)} &= 
\prod _{n=0}^{\infty }{\frac {(1-q^{8n+3})(1-q^{8n+5})}{(1-q^{8n+1})(1-q^{8n+7})}}\\
&=1 + {\cfrac
{q+q^2}{1+{\cfrac {q^4}{1+{\cfrac {q^{3}+q^6}{1+{\cfrac
{q^{8}}{1+\ddots }}}}}}}}.
\end{align*}               

Some of the identities given in this section are due Robins \cite{Ro-thesis},
and many are due to Huang \cite{Hu1998}.  Any identities that appear
to be new are marked $\star$.

\subsubsection{Type 1}
\mylabel{subsubsec:type8.1}
\begin{alignat}{3}
G(3)\,H(1) - G(1)\,H(3) &= 
\frac{\eta(\tau)\eta(12\tau)^{2}}{\eta(3\tau)\eta(8\tau)\eta(24\tau)},&\quad &\Gamma_1(24), &\quad -B=6,
\mylabel{eq:RR811}\\
G(3)\,H(1) + G(1)\,H(3) &= 
\frac{\eta(2\tau)\eta(4\tau)^{2}\eta(6\tau)^{2}}{\eta(\tau)\eta(3\tau)\eta(8\tau)^{2}\eta(12\tau)},&\quad &\Gamma_1(24), &\quad -B=6,
\mylabel{eq:RR812}\\
G(5)\,H(1) - G(1)\,H(5) &= 
\frac{\eta(2\tau)\eta(10\tau)\eta(20\tau)}{\eta(5\tau)\eta(8\tau)\eta(40\tau)},&\quad &\Gamma_1(40), &\quad -B=20,
\mylabel{eq:RR813}\\
G(7)\,H(1) - G(1)\,H(7) &= 
\frac{\eta(4\tau)\eta(28\tau)}{\eta(8\tau)\eta(56\tau)},&\quad &\Gamma_1(56), &\quad -B=36,
\mylabel{eq:RR814}\\
G(9)\,H(1) - G(1)\,H(9) &= 
\frac{\eta(4\tau)\eta(6\tau)^{2}\eta(36\tau)}{\eta(3\tau)\eta(8\tau)\eta(12\tau)\eta(72\tau)},&\quad &\Gamma_1(72), &\quad -B=60,
\mylabel{eq:RR815}\\
G(5)\,H(3) - G(3)\,H(5) &= 
\frac{\eta(\tau)\eta(4\tau)\eta(6\tau)\eta(10\tau)\eta(15\tau)\eta(60\tau)}{\eta(2\tau)\eta(3\tau)\eta(5\tau)\eta(24\tau)\eta(30\tau)\eta(40\tau)},&\quad &\Gamma_1(120), &\quad -B=144.
\mylabel{eq:RR816}
\end{alignat}

\subsubsection{Type 2}
\mylabel{subsubsec:type8.2}
\begin{alignat}{3}
G(1)\,G(1) - H(1)\,H(1) &= 
\frac{\eta(4\tau)^{6}}{\eta(\tau)\eta(2\tau)\eta(8\tau)^{4}},&\quad &\Gamma_1(8), &\quad -B=1,
\mylabel{eq:RR821}\\
G(1)\,G(1) + H(1)\,H(1) &= 
\frac{\eta(2\tau)^{6}}{\eta(\tau)^{3}\eta(4\tau)\eta(8\tau)^{2}},&\quad &\Gamma_1(8), &\quad -B=1,
\mylabel{eq:RR822}\\
G(1)\,G(3) - H(1)\,H(3) &= 
\frac{\eta(2\tau)^{2}\eta(6\tau)\eta(12\tau)^{2}}{\eta(\tau)\eta(3\tau)\eta(4\tau)\eta(24\tau)^{2}},&\quad &\Gamma_1(24), &\quad -B=6,
\mylabel{eq:RR823}\\
G(1)\,G(3) + H(1)\,H(3) &= 
\frac{\eta(3\tau)\eta(4\tau)^{2}}{\eta(\tau)\eta(8\tau)\eta(24\tau)},&\quad &\Gamma_1(24), &\quad -B=6,
\mylabel{eq:RR824}\\
G(1)\,G(5) + H(1)\,H(5) &= 
\frac{\eta(2\tau)\eta(4\tau)\eta(10\tau)}{\eta(\tau)\eta(8\tau)\eta(40\tau)},&\quad &\Gamma_1(40), &\quad -B=20,
\mylabel{eq:RR825}\\
G(1)\,G(9) + H(1)\,H(9) &= 
\frac{\eta(2\tau)\eta(3\tau)\eta(12\tau)\eta(18\tau)}{\eta(\tau)\eta(8\tau)\eta(9\tau)\eta(72\tau)},&\quad &\Gamma_1(72), &\quad -B=60,
\mylabel{eq:RR826}\\
G(1)\,G(15) + H(1)\,H(15) &= 
\frac{\eta(2\tau)\eta(3\tau)\eta(5\tau)\eta(12\tau)\eta(20\tau)\eta(30\tau)}{\eta(\tau)\eta(6\tau)\eta(8\tau)\eta(10\tau)\eta(15\tau)\eta(120\tau)},&\quad &\Gamma_1(120), &\quad -B=144.
\mylabel{eq:RR827}
\end{alignat}

\subsubsection{Type 3}
\mylabel{subsubsec:type8.3}
\begin{alignat}{3}
\star\quad\frac{G(3)\,G(5) - H(3)\,H(5)}{G(15)\,H(1) + H(15)\,G(1)} &= 
\frac{\eta(4\tau)\eta(60\tau)}{\eta(12\tau)\eta(20\tau)},&\quad &\Gamma_1(120), &\quad -B=256,
\mylabel{eq:RR831}\\
\frac{G(3)\,G(5) + H(3)\,H(5)}{G(15)\,H(1) - H(15)\,G(1)} &= 
\frac{\eta(8\tau)\eta(12\tau)\eta(20\tau)\eta(120\tau)}{\eta(4\tau)\eta(24\tau)\eta(40\tau)\eta(60\tau)},&\quad &\Gamma_1(120), &\quad -B=224,
\mylabel{eq:RR832}\\
\star\quad\frac{G(1)\,G(15) - H(1)\,H(15)}{G(5)\,H(3) + H(5)\,G(3)} &= 
\frac{\eta(4\tau)^{2}\eta(6\tau)\eta(10\tau)\eta(24\tau)^{2}\eta(40\tau)^{2}\eta(60\tau)^{2}}{\eta(2\tau)\eta(8\tau)^{2}\eta(12\tau)^{2}\eta(20\tau)^{2}\eta(30\tau)\eta(120\tau)^{2}},&\quad &\Gamma_1(120), &\quad -B=192,
\mylabel{eq:RR833}\\
\star\quad\frac{G(3)\,G(7) + H(3)\,H(7)}{G(21)\,H(1) - H(21)\,G(1)} &= 
\frac{\eta(8\tau)\eta(21\tau)\eta(28\tau)\eta(168\tau)}{\eta(7\tau)\eta(24\tau)\eta(56\tau)\eta(84\tau)},&\quad &\Gamma_1(168), &\quad -B=528,
\mylabel{eq:RR834}\\
\star\quad\frac{G(1)\,G(21) + H(1)\,H(21)}{G(7)\,H(3) - H(7)\,G(3)} &= 
\frac{\eta(3\tau)\eta(4\tau)\eta(24\tau)\eta(56\tau)}{\eta(\tau)\eta(8\tau)\eta(12\tau)\eta(168\tau)},&\quad &\Gamma_1(168), &\quad -B=528,
\mylabel{eq:RR835}\\
\frac{G(1)\,G(39) + H(1)\,H(39)}{G(13)\,H(3) - H(13)\,G(3)} &= 
\frac{\eta(2\tau)\eta(3\tau)\eta(13\tau)\eta(24\tau)\eta(78\tau)\eta(104\tau)}{\eta(\tau)\eta(6\tau)\eta(8\tau)\eta(26\tau)\eta(39\tau)\eta(312\tau)},&\quad &\Gamma_1(312), &\quad -B=1632,
\mylabel{eq:RR836}\\
\frac{G(1)\,G(55) + H(1)\,H(55)}{G(11)\,H(5) - H(11)\,G(5)} &= 
\frac{\eta(2\tau)\eta(5\tau)\eta(11\tau)\eta(40\tau)\eta(88\tau)\eta(110\tau)}{\eta(\tau)\eta(8\tau)\eta(10\tau)\eta(22\tau)\eta(55\tau)\eta(440\tau)},&\quad &\Gamma_1(440), &\quad -B=3680.
\mylabel{eq:RR837}
\end{alignat}

\subsubsection{Type 8}
\mylabel{subsubsec:type8.8}
\begin{equation}
\star\quad G(1)^{3}\,H(3) - H(1)^{3}\,G(3) = 
{3}\,\frac{\eta(2\tau)^{3}\eta(4\tau)\eta(6\tau)\eta(24\tau)^{2}}{\eta(\tau)^{2}\eta(3\tau)\eta(8\tau)^{4}},\quad \Gamma_1(24), \quad -B=10.
\mylabel{eq:RR881}
\end{equation}


\subsection{Mod $10$}
\mylabel{subsec:mod10}
There is only one real non-principal character mod $10$ that satisfies 
$\chi(-1)=1$,
namely the character $\chi_{10}$ induced by the Legendre symbol mod $5$.
In this section
$$
G(1) = G(1,10,\chi_{10}) 
= \frac{1}{\eta_{10;1}(\tau)} = \frac{q^{-23/60}}{(q,q^9;q^{10})_\infty},
\quad
H(1)= H(1,10,\chi_{10})
= \frac{1}{\eta_{10;3}(\tau)} = \frac{q^{13/60}}{(q^3,q^{7};q^{10})_\infty}.
$$
All the identities in this section appear to be new.

\subsubsection{Type 1}
\mylabel{subsubsec:type10.1}
\begin{alignat}{3}
G(6)\,H(1) - G(1)\,H(6) &= 
\frac{\eta(4\tau)\eta(5\tau)\eta(12\tau)\eta(30\tau)^{3}}{\eta(6\tau)\eta(10\tau)^{2}\eta(15\tau)\eta(60\tau)^{2}},&\quad &\Gamma_1(60), &\quad -B=40.
\mylabel{eq:RR1011}
\end{alignat}

\subsubsection{Type 2}
\mylabel{subsubsec:type10.2}
\begin{alignat}{3}
G(2)\,G(3) - H(2)\,H(3) &= 
\frac{\eta(4\tau)\eta(10\tau)^{3}\eta(12\tau)\eta(15\tau)}{\eta(2\tau)\eta(5\tau)\eta(20\tau)^{2}\eta(30\tau)^{2}},&\quad &\Gamma_1(60), &\quad -B=40,
\mylabel{eq:RR1021}\\
G(1)\,G(9) - H(1)\,H(9) &= 
\frac{\eta(2\tau)\eta(3\tau)\eta(5\tau)\eta(18\tau)\eta(30\tau)^{2}\eta(45\tau)}{\eta(\tau)\eta(9\tau)\eta(10\tau)^{2}\eta(15\tau)\eta(90\tau)^{2}},&\quad &\Gamma_1(90), &\quad -B=96.
\mylabel{eq:RR1022}
\end{alignat}

\subsubsection{Type 5}
\mylabel{subsubsec:type10.5}
\begin{equation}
G^{*}(1)\,G^{*}(4) - H^{*}(1)\,H^{*}(4) = 
\frac{\eta(\tau)\eta(4\tau)^{3}\eta(10\tau)\eta(16\tau)\eta(40\tau)}{\eta(2\tau)^{2}\eta(5\tau)\eta(8\tau)^{2}\eta(20\tau)\eta(80\tau)},\quad \Gamma_1(80), \quad -B=64.
\mylabel{eq:RR1051}
\end{equation}

\subsubsection{Type 6}
\mylabel{subsubsec:type10.6}
\begin{alignat}{3}
G(1)\,H^{*}(1) - G^{*}(1)\,H(1) &= 
{2}\,\frac{\eta(20\tau)^{2}}{\eta(10\tau)^{2}},&\quad &\Gamma_1(20), &\quad -B=4,
\mylabel{eq:RR1061}\\
G(1)\,H^{*}(1) + G^{*}(1)\,H(1) &= 
{2}\,\frac{\eta(4\tau)^{2}}{\eta(2\tau)\eta(10\tau)},&\quad &\Gamma_1(20), &\quad -B=4.
\mylabel{eq:RR1062}
\end{alignat}

\subsubsection{Type 8}
\mylabel{subsubsec:type10.8}
\begin{alignat}{3}
G(1)^{2}\,H(2) - H(1)^{2}\,G(2) = 
{2}\,\frac{\eta(2\tau)\eta(5\tau)\eta(20\tau)^{2}}{\eta(\tau)\eta(10\tau)^{3}},\quad \Gamma_1(20), \quad -B=4,
\mylabel{eq:RR1081}\\
G(1)^{2}\,H(2) + H(1)^{2}\,G(2) = 
{2}\,\frac{\eta(4\tau)^{2}\eta(5\tau)}{\eta(\tau)\eta(10\tau)^{2}},\quad \Gamma_1(20), \quad -B=4,
\mylabel{eq:RR1082}\\
G(1)^{3}\,H(3) - H(1)^{3}\,G(3) = 
{3}\,\frac{\eta(2\tau)^{3}\eta(5\tau)^{2}\eta(6\tau)\eta(15\tau)\eta(30\tau)}{\eta(\tau)^{2}\eta(3\tau)\eta(10\tau)^{5}},\quad \Gamma_1(30), \quad -B=16.
\mylabel{eq:RR1083}
\end{alignat}

\subsection{Mod $12$}
\mylabel{subsec:mod12}
There is only one non-principal character mod $12$ that satisfies 
$\chi(-1)=1$, namely $\chi(\cdot)= \leg{12}{\cdot}$.
In this section
$$
G(1)=G(1,12,\leg{12}{\cdot}) = \frac{1}{\eta_{12;1}(\tau)} 
= \frac{q^{-13/24}}{(q,q^{11};q^{12})_\infty},
\quad
H(1)=H(1,12,\leg{12}{\cdot}) = \frac{1}{\eta_{12;5}(\tau)} 
= \frac{q^{11/24}}{(q^5,q^7;q^{12})_\infty}.
$$
These functions were considered by Robins \cite[p17]{Ro-thesis},
who found \eqn{RR1213}, \eqn{RR1214}, \eqn{RR1221}, \eqn{RR1222}.
The remaining identities appear to be new and are marked $\star$.

\subsubsection{Type 1}
\mylabel{subsubsec:type12.1}
\begin{alignat}{3}
\star\quad G(2)\,H(1) - G(1)\,H(2) &= 
\frac{\eta(\tau)\eta(4\tau)\eta(6\tau)}{\eta(2\tau)\eta(12\tau)^{2}},&\quad &\Gamma_1(24), &\quad -B=4,
\mylabel{eq:RR1211}\\
\star\quad G(2)\,H(1) + G(1)\,H(2) &= 
\frac{\eta(3\tau)^{2}\eta(4\tau)}{\eta(\tau)\eta(12\tau)^{2}},&\quad &\Gamma_1(24), &\quad -B=4,
\mylabel{eq:RR1212}\\
G(3)\,H(1) - G(1)\,H(3) &= 
\frac{\eta(2\tau)\eta(18\tau)}{\eta(12\tau)\eta(36\tau)},&\quad &\Gamma_1(36), &\quad -B=12,
\mylabel{eq:RR1213}\\
G(3)\,H(1) + G(1)\,H(3) &= 
\frac{\eta(4\tau)\eta(6\tau)^{5}\eta(9\tau)^{2}}{\eta(2\tau)\eta(3\tau)^{2}\eta(12\tau)^{3}\eta(18\tau)^{2}},&\quad &\Gamma_1(36), &\quad -B=12,
\mylabel{eq:RR1214}\\
\star\quad G(4)\,H(1) - G(1)\,H(4) &= 
\frac{\eta(3\tau)\eta(16\tau)}{\eta(12\tau)\eta(48\tau)},&\quad &\Gamma_1(48), &\quad -B=24,
\mylabel{eq:RR1215}\\
\star\quad G(5)\,H(1) - G(1)\,H(5) &= 
\frac{\eta(4\tau)\eta(6\tau)\eta(10\tau)\eta(15\tau)}{\eta(5\tau)\eta(12\tau)^{2}\eta(60\tau)},&\quad &\Gamma_1(60), &\quad -B=40,
\mylabel{eq:RR1216}\\
\star\quad G(3)\,H(2) - G(2)\,H(3) &= 
\frac{\eta(\tau)\eta(6\tau)\eta(8\tau)\eta(9\tau)\eta(12\tau)\eta(72\tau)}{\eta(2\tau)\eta(3\tau)\eta(24\tau)^{2}\eta(36\tau)^{2}},&\quad &\Gamma_1(72), &\quad -B=48,
\mylabel{eq:RR1217}\\
\star\quad G(6)\,H(1) - G(1)\,H(6) &= 
\frac{\eta(8\tau)\eta(9\tau)}{\eta(12\tau)\eta(72\tau)},&\quad &\Gamma_1(72), &\quad -B=60.
\mylabel{eq:RR1218}
\end{alignat}

\subsubsection{Type 2}
\mylabel{subsubsec:type12.2}
\begin{alignat}{3}
G(1)^2 - H(1)^2 &= 
\frac{\eta(2\tau)^{3}\eta(6\tau)^{3}}{\eta(\tau)^{2}\eta(12\tau)^{4}},&\quad &\Gamma_1(12), &\quad -B=2,
\mylabel{eq:RR1221}\\
G(1)^2 + H(1)^2 &= 
\frac{\eta(2\tau)\eta(3\tau)^{4}\eta(4\tau)}{\eta(\tau)^{2}\eta(6\tau)\eta(12\tau)^{3}},&\quad &\Gamma_1(12), &\quad -B=2,
\mylabel{eq:RR1222}\\
\star\quad G(1)\,G(2) - H(1)\,H(2) &= 
\frac{\eta(3\tau)^{2}\eta(8\tau)^{2}}{\eta(\tau)\eta(12\tau)\eta(24\tau)^{2}},&\quad &\Gamma_1(24), &\quad -B=8,
\mylabel{eq:RR1223}\\
\star\quad G(1)\,G(3) - H(1)\,H(3) &= 
\frac{\eta(2\tau)\eta(4\tau)\eta(9\tau)\eta(18\tau)}{\eta(\tau)\eta(12\tau)\eta(36\tau)^{2}},&\quad &\Gamma_1(36), &\quad -B=18,
\mylabel{eq:RR1224}\\
\star\quad G(1)\,G(5) - H(1)\,H(5) &= 
\frac{\eta(2\tau)\eta(3\tau)\eta(20\tau)\eta(30\tau)}{\eta(\tau)\eta(12\tau)\eta(60\tau)^{2}},&\quad &\Gamma_1(60), &\quad -B=40.
\mylabel{eq:RR1225}
\end{alignat}

\subsubsection{Type 3}
\mylabel{subsubsec:type12.3}
\begin{alignat}{3}
\frac{G(1)\,G(10) - H(1)\,H(10)}{G(5)\,H(2) - H(5)\,G(2)} &= 
\frac{\eta(2\tau)\eta(5\tau)\eta(24\tau)^{2}\eta(60\tau)^{2}}{\eta(\tau)\eta(10\tau)\eta(12\tau)^{2}\eta(120\tau)^{2}},&\quad &\Gamma_1(120), &\quad -B=256,
\mylabel{eq:RR1231}\\
\frac{G(5)\,G(7) - H(5)\,H(7)}{G(35)\,H(1) - H(35)\,G(1)} &= 
\frac{\eta(12\tau)\eta(420\tau)}{\eta(60\tau)\eta(84\tau)},&\quad &\Gamma_1(420), &\quad -B=3648,
\mylabel{eq:RR1232}\\
\frac{G(1)\,G(35) - H(1)\,H(35)}{G(7)\,H(5) - H(7)\,G(5)} &= 
\frac{\eta(3\tau)\eta(4\tau)\eta(5\tau)\eta(7\tau)\eta(60\tau)^{2}\eta(84\tau)^{2}\eta(105\tau)\eta(140\tau)}{\eta(\tau)\eta(12\tau)^{2}\eta(15\tau)\eta(20\tau)\eta(21\tau)\eta(28\tau)\eta(35\tau)\eta(420\tau)^{2}},&\quad &\Gamma_1(420), &\quad -B=2880.
\mylabel{eq:RR1233}
\end{alignat}

\subsubsection{Type 4}
\mylabel{subsubsec:type12.4}
\begin{equation}
\star\quad G^{*}(2)\,H^{*}(1) - G^{*}(1)\,H^{*}(2) = 
\frac{\eta(\tau)\eta(6\tau)\eta(8\tau)^{3}\eta(48\tau)}{\eta(2\tau)\eta(4\tau)\eta(16\tau)\eta(24\tau)^{3}},\quad \Gamma_1(48), \quad -B=24.
\mylabel{eq:RR1241}
\end{equation}

\subsubsection{Type 5}
\mylabel{subsubsec:type12.5}
\begin{equation}
\star\quad G^{*}(1)\,G^{*}(2) - H^{*}(1)\,H^{*}(2) = 
\frac{\eta(\tau)\eta(6\tau)\eta(16\tau)}{\eta(2\tau)\eta(12\tau)\eta(48\tau)},\quad \Gamma_1(48), \quad -B=24.
\mylabel{eq:RR1251}
\end{equation}

\subsubsection{Type 6}
\mylabel{subsubsec:type12.6}
\begin{alignat}{3}
\star\quad G(1)\,H^{*}(1) - G^{*}(1)\,H(1) &= 
{2}\,\frac{\eta(4\tau)^{2}\eta(6\tau)^{2}\eta(24\tau)^{3}}{\eta(2\tau)\eta(8\tau)\eta(12\tau)^{5}},&\quad &\Gamma_1(24), &\quad -B=4,
\mylabel{eq:RR1261}\\
\star\quad G(1)\,H^{*}(1) + G^{*}(1)\,H(1) &= 
{2}\,\frac{\eta(4\tau)\eta(6\tau)^{2}\eta(8\tau)\eta(24\tau)}{\eta(2\tau)\eta(12\tau)^{4}},&\quad &\Gamma_1(24), &\quad -B=4.
\mylabel{eq:RR1262}
\end{alignat}

\subsubsection{Type 7}
\mylabel{subsubsec:type12.7}
\begin{alignat}{3}
\star\quad G^{*}(1)\,G(1) - H^{*}(1)\,H(1) = 
\frac{\eta(8\tau)\eta(12\tau)^{4}}{\eta(4\tau)\eta(6\tau)\eta(24\tau)^{3}},\quad \Gamma_1(24), \quad -B=4,
\mylabel{eq:RR1271}\\
\star\quad G^{*}(1)\,G(1) + H^{*}(1)\,H(1) = 
\frac{\eta(4\tau)^{4}\eta(6\tau)}{\eta(2\tau)^{2}\eta(8\tau)\eta(12\tau)\eta(24\tau)},\quad \Gamma_1(24), \quad -B=4.
\mylabel{eq:RR1272}
\end{alignat}

\subsubsection{Type 8}
\mylabel{subsubsec:type12.8}
\begin{alignat}{3}
\star\quad G(1)^{2}\,H(2) - H(1)^{2}\,G(2) = 
{2}\,\frac{\eta(3\tau)\eta(4\tau)^{2}\eta(6\tau)\eta(24\tau)^{3}}{\eta(\tau)\eta(8\tau)\eta(12\tau)^{5}},\quad \Gamma_1(24), \quad -B=4,
\mylabel{eq:RR1281}\\
\star\quad G(1)^{2}\,H(2) + H(1)^{2}\,G(2) = 
{2}\,\frac{\eta(3\tau)\eta(4\tau)\eta(6\tau)\eta(8\tau)\eta(24\tau)}{\eta(\tau)\eta(12\tau)^{4}},\quad \Gamma_1(24), \quad -B=4,
\mylabel{eq:RR1282}\\
\star\quad G(1)^{3}\,H(3) - H(1)^{3}\,G(3) = 
{3}\,\frac{\eta(2\tau)^{2}\eta(3\tau)\eta(4\tau)\eta(6\tau)\eta(9\tau)\eta(36\tau)^{2}}{\eta(\tau)^{2}\eta(12\tau)^{5}\eta(18\tau)},\quad \Gamma_1(36), \quad -B=18.
\mylabel{eq:RR1283}
\end{alignat}


\subsection{Mod $13$}
There is only one non-principal character mod $13$ that satisfies $\chi(-1)=1$,
namely $\chi(\cdot)= \leg{\cdot}{13}$.
In this section
\begin{align*}
G(1)&=G(1,13,\leg{\cdot}{13})= \frac{1}
{
\eta_{13; 1,3,4}(\tau)
}
=
\frac{q^{-1/4}}
{
(q,q^{3},q^{4},q^{9},q^{10},q^{12};q^{13})_\infty
},\\
H(1)&=H(1,13,\leg{\cdot}{13})= \frac{1}
{
\eta_{13; 2,5,6}(\tau)
}
=
\frac{q^{3/4}}
{
(q^{2},q^{5},q^{6},q^{7},q^{8},q^{11};q^{13})_\infty
}.
\end{align*}
These functions were considered by Robins \cite[p.18]{Ro-thesis},
who found the one identity \eqn{RR1311}. The remaining four identities
appear to be new and are marked $\star$.

\subsubsection{Type 1}
\mylabel{subsubsec:type13.1}
\begin{alignat}{3}
G(3)\,H(1) - G(1)\,H(3) &= 
1
,&\quad &\Gamma_1(39), &\quad -B=24.
\mylabel{eq:RR1311}
\end{alignat}

\subsubsection{Type 3}
\mylabel{subsubsec:type13.3}
\begin{alignat}{3}
\star\quad \frac{G(1)\,G(2) + H(1)\,H(2)}{G(2)\,H(1) - H(2)\,G(1)} &= 
\frac{\eta(2\tau)^{2}\eta(13\tau)^{2}}{\eta(\tau)^{2}\eta(26\tau)^{2}},&\quad &\Gamma_1(26), &\quad -B=18,
\mylabel{eq:RR1331}\\
\star\quad \frac{G(2)\,G(5) + H(2)\,H(5)}{G(10)\,H(1) - H(10)\,G(1)} &= 
1
,&\quad &\Gamma_1(130), &\quad -B=432,
\mylabel{eq:RR1332}\\
\star\quad \frac{G(1)\,G(14) + H(1)\,H(14)}{G(7)\,H(2) - H(7)\,G(2)} &= 
\frac{\eta(2\tau)\eta(7\tau)\eta(26\tau)\eta(91\tau)}{\eta(\tau)\eta(13\tau)\eta(14\tau)\eta(182\tau)},&\quad &\Gamma_1(182), &\quad -B=864.
\mylabel{eq:RR1333}
\end{alignat}

\subsubsection{Type 9}
\mylabel{subsubsec:type13.9}
\begin{equation}
\star\quad G(1)^{3}\,H(1) - H(1)^{3}\,G(1) - 1 = 
{3}\,\frac{\eta(13\tau)^{2}}{\eta(\tau)^{2}},\quad \Gamma_1(13), \quad -B=6.
\mylabel{eq:RR1391}
\end{equation}

\subsection{Mod $15$}
There is only one real non-principal character mod $15$ that satisfies 
$\chi(-1)=1$,
namely the one induced by the Legendre symbol mod $5$:
$$
\chi_{15}(n) =
\begin{cases}
1, & n\equiv \pm1,4 \pmod{15}, \\
-1, & n\equiv \pm2,7 \pmod{15}, \\
0, & \mbox{otherwise}.                    
\end{cases}
$$
Thus in this section
\begin{align*}
G(1)&=G(1,15,\chi_{15})= \frac{1}
{
\eta_{15; 1,4}(\tau)
}
=
\frac{q^{-17/30}}
{
(q,q^{4},q^{11},q^{14};q^{15})_\infty
},\\
H(1)&=H(1,15,\chi_{15})= \frac{1}
{
\eta_{15; 2,7}(\tau)
}
=
\frac{q^{7/30}}
{
(q^{2},q^{7},q^{8},q^{13};q^{15})_\infty
}.
\end{align*}
All the identities in this section appear to be new.

\subsubsection{Type 2}
\mylabel{subsubsec:type15.2}
\begin{alignat}{3}
G(1)\,G(4) - H(1)\,H(4) &= 
\frac{\eta(2\tau)\eta(3\tau)\eta(10\tau)\eta(12\tau)\eta(30\tau)^{2}}{\eta(\tau)\eta(4\tau)\eta(15\tau)^{2}\eta(60\tau)^{2}},&\quad &\Gamma_1(60), &\quad -B=48.
\mylabel{eq:RR1521}
\end{alignat}

\subsubsection{Type 3}
\mylabel{subsubsec:type15.3}
\begin{alignat}{3}
\frac{G(2)\,G(3) - H(2)\,H(3)}{G(6)\,H(1) - H(6)\,G(1)} &= 
\frac{\eta(6\tau)\eta(10\tau)\eta(15\tau)^{3}\eta(90\tau)}{\eta(3\tau)\eta(5\tau)\eta(30\tau)^{3}\eta(45\tau)},&\quad &\Gamma_1(90), &\quad -B=120.
\mylabel{eq:RR1531}
\end{alignat}

\subsubsection{Type 6}
\mylabel{subsubsec:type15.6}
\begin{equation}
G(1)\,H^{*}(1) - G^{*}(1)\,H(1) = 
{2}\,\frac{\eta(4\tau)\eta(6\tau)^{3}\eta(10\tau)\eta(60\tau)^{2}}{\eta(2\tau)^{2}\eta(12\tau)\eta(30\tau)^{4}},\quad \Gamma_1(60), \quad -B=48.
\mylabel{eq:RR1561}
\end{equation}

\subsubsection{Type 8}
\mylabel{subsubsec:type15.8}
\begin{equation}
G(1)^{2}\,H(2) + H(1)^{2}\,G(2) = 
{2}\,\frac{\eta(3\tau)^{2}\eta(6\tau)\eta(10\tau)^{2}}{\eta(\tau)\eta(2\tau)\eta(15\tau)^{3}},\quad \Gamma_1(30), \quad -B=12.
\mylabel{eq:RR1581}
\end{equation}

\subsection{Mod $17$}
There is only one non-principal character mod $17$ that satisfies $\chi(-1)=1$,
namely $\chi(\cdot)= \leg{\cdot}{17}$. In this section
\begin{align*}
G(1)&=G(1,17,\leg{\cdot}{17})= \frac{1}
{
\eta_{17; 1,2,4,8}(\tau)
}
\\
&=
\frac{q^{-2/3}}
{
(q,q^{2},q^{4},q^{8},q^{9},q^{13},q^{15},q^{16};q^{17})_\infty
},
\end{align*}
\begin{align*}
H(1)&=H(1,17,\leg{\cdot}{17})= \frac{1}
{
\eta_{17; 3,5,6,7}(\tau)
}
\\
&=
\frac{q^{4/3}}
{
(q^{3},q^{5},q^{6},q^{7},q^{10},q^{11},q^{12},q^{14};q^{17})_\infty
}.
\end{align*}
These functions were not considered by Robins \cite{Ro-thesis}.
Nonetheless we find one identity.

\subsubsection{Type 1}
\mylabel{subsubsec:type17.1}
\begin{alignat}{3}
G(2)\,H(1) - G(1)\,H(2) &= 
1
,&\quad &\Gamma_1(34), &\quad -B=16.
\mylabel{eq:RR1711}
\end{alignat}

\subsection{Mod $21$}
\mylabel{subsec:mod21}
There is only one non-principal character mod $21$ that satisfies $\chi(-1)=1$,
namely $\chi(\cdot)= \leg{21}{\cdot}$.
In this section 
$$                   
G(1) = G(1,21,\leg{21}{\cdot}) = 
\frac{1}{
\eta_{21; 1,4,5}(\tau)
}
=
\frac{q^{-5/4}}
{
(q,q^{4},q^{5},q^{16},q^{17},q^{20};q^{21})_\infty
},
$$                 
$$                   
H(1) = H(1,21,\leg{21}{\cdot}) = 
\frac{1}{
\eta_{21; 2,8,10}(\tau)
}
=
\frac{q^{3/4}}
{
 (q^{2},q^{8},q^{10},q^{11},q^{13},q^{19};q^{21})_\infty
}.
$$

\subsubsection{Type 1}
\mylabel{subsubsec:type21.1}
\begin{alignat}{3}
G(2)\,H(1) - G(1)\,H(2) &= 
\frac{\eta(3\tau)\eta(6\tau)\eta(7\tau)^{2}}{\eta(2\tau)\eta(21\tau)^{3}},&\quad &\Gamma_1(42), &\quad -B=24,
\mylabel{eq:RR2111}\\
G(4)\,H(1) - G(1)\,H(4) &= 
\frac{\eta(6\tau)^{2}\eta(7\tau)\eta(28\tau)}{\eta(2\tau)\eta(21\tau)\eta(42\tau)\eta(84\tau)},&\quad &\Gamma_1(84), &\quad -B=96.
\mylabel{eq:RR2112}
\end{alignat}

\subsubsection{Type 2}
\mylabel{subsubsec:type21.2}
\begin{alignat}{3}
G(1)\,G(2) - H(1)\,H(2) &= 
\frac{\eta(3\tau)\eta(6\tau)\eta(14\tau)^{2}}{\eta(\tau)\eta(42\tau)^{3}},&\quad &\Gamma_1(42), &\quad -B=24.
\mylabel{eq:RR2121}
\end{alignat}

\subsubsection{Type 7}
\mylabel{subsubsec:type21.7}
\begin{equation}
G^{*}(1)\,G(1) - H^{*}(1)\,H(1) = 
\frac{\eta(6\tau)^{2}\eta(14\tau)^{3}\eta(84\tau)}{\eta(2\tau)\eta(28\tau)\eta(42\tau)^{4}},\quad \Gamma_1(84), \quad -B=96.
\mylabel{eq:RR2171}
\end{equation}

\subsection{Mod $24$}
There are three real non-principal characters mod $24$ that satisfy   
$\chi(-1)=1$.

\begin{enumerate}
\item[(i)]
The character $\chi_{24,1}(\cdot)$ induced by $\leg{8}{\cdot}$.
\item[(ii)]
The character $\chi_{24,2}(\cdot)=\leg{12}{\cdot}$ covered previously
in Section \subsect{mod12}.
\item[(iii)]
The character $\chi_{24,3}(\cdot)=\leg{24}{\cdot}$.
\end{enumerate}

\subsubsection{$\chi_{24,1}$}
\mylabel{subsubsec:chi241}
We have 
\begin{align*}
G(1)&=G(1,24,\chi_{24,1})= \frac{1}
{
\eta_{24; 1,7}(\tau)
}
=
\frac{q^{-25/24}}
{
(q,q^{7},q^{17},q^{23};q^{24})_\infty
},\\
H(1)&=H(1,24,\chi_{24,1})= \frac{1}
{
\eta_{24; 5,11}(\tau)
}
=
\frac{q^{23/24}}
{
(q^{5},q^{11},q^{13},q^{19};q^{24})_\infty
}.
\end{align*}

\paragraph[0pt]{Type 1}
\mylabel{subsubsubsec:type24a.1}
\begin{alignat}{3}
G(2)\,H(1) - G(1)\,H(2) &= 
\frac{\eta(3\tau)\eta(12\tau)^{2}}{\eta(6\tau)\eta(24\tau)^{2}},&\quad &\Gamma_1(48), &\quad -B=24,
\mylabel{eq:RR24a11}\\
G(2)\,H(1) + G(1)\,H(2) &= 
\frac{\eta(4\tau)^{3}\eta(6\tau)^{4}}{\eta(2\tau)^{2}\eta(3\tau)\eta(8\tau)\eta(12\tau)^{2}\eta(24\tau)},&\quad &\Gamma_1(48), &\quad -B=24,
\mylabel{eq:RR24a12}\\
G(3)\,H(1) - G(1)\,H(3) &= 
\frac{\eta(4\tau)\eta(6\tau)^{2}\eta(9\tau)\eta(36\tau)}{\eta(3\tau)\eta(8\tau)\eta(12\tau)\eta(18\tau)\eta(72\tau)},&\quad &\Gamma_1(72), &\quad -B=60.
\mylabel{eq:RR24a13}
\end{alignat}

\paragraph[0pt]{Type 2}
\mylabel{subsubsubsec:type24a.2}
\begin{alignat}{3}
G(1)\,G(1) - H(1)\,H(1) &= 
\frac{\eta(2\tau)^{2}\eta(3\tau)^{2}\eta(4\tau)\eta(12\tau)^{2}}{\eta(\tau)^{2}\eta(6\tau)\eta(8\tau)\eta(24\tau)^{3}},&\quad &\Gamma_1(24), &\quad -B=12.
\mylabel{eq:RR24a21}
\end{alignat}

\paragraph[0pt]{Type 6}
\mylabel{subsubsubsec:type24a.6}
\begin{alignat}{3}
G(1)\,H^{*}(1) - G^{*}(1)\,H(1) &= 
{2}\,\frac{\eta(4\tau)^{2}\eta(12\tau)^{2}\eta(48\tau)^{2}}{\eta(2\tau)\eta(8\tau)\eta(24\tau)^{4}},&\quad &\Gamma_1(48), &\quad -B=24,
\mylabel{eq:RR24a61}\\
G(1)\,H^{*}(1) + G^{*}(1)\,H(1) &= 
{2}\,\frac{\eta(4\tau)\eta(6\tau)^{2}\eta(16\tau)\eta(48\tau)}{\eta(2\tau)\eta(8\tau)\eta(24\tau)^{3}},&\quad &\Gamma_1(48), &\quad -B=24.
\mylabel{eq:RR24a62}
\end{alignat}

\paragraph[0pt]{Type 7}
\mylabel{subsubsubsec:type24a.7}
\begin{alignat}{3}
G^{*}(1)\,G(1) - H^{*}(1)\,H(1) = 
\frac{\eta(4\tau)^{2}\eta(24\tau)^{2}}{\eta(2\tau)\eta(8\tau)\eta(48\tau)^{2}},\quad \Gamma_1(48), \quad -B=24,
\mylabel{eq:RR24a71}\\
G^{*}(1)\,G(1) + H^{*}(1)\,H(1) = 
\frac{\eta(6\tau)^{2}\eta(8\tau)^{2}}{\eta(2\tau)\eta(12\tau)\eta(16\tau)\eta(48\tau)},\quad \Gamma_1(48), \quad -B=24.
\mylabel{eq:RR24a72}
\end{alignat}

\paragraph[0pt]{Type 8}
\mylabel{subsubsubsec:type24a.8}
\begin{alignat}{3}
G(1)^{2}\,H(2) - H(1)^{2}\,G(2) = 
{2}\,\frac{\eta(3\tau)\eta(4\tau)^{2}\eta(12\tau)^{2}\eta(48\tau)^{2}}{\eta(\tau)\eta(6\tau)\eta(8\tau)\eta(24\tau)^{4}},\quad \Gamma_1(48), \quad -B=24,
\mylabel{eq:RR24a81}\\
G(1)^{2}\,H(2) + H(1)^{2}\,G(2) = 
{2}\,\frac{\eta(3\tau)\eta(4\tau)\eta(6\tau)\eta(16\tau)\eta(48\tau)}{\eta(\tau)\eta(8\tau)\eta(24\tau)^{3}},\quad \Gamma_1(48), \quad -B=24.
\mylabel{eq:RR24a82}
\end{alignat}

\paragraph[0pt]{Type 10}
\mylabel{subsubsubsec:type24a.10}
\begin{equation}
\frac{G(3)\,H(2) + H(3)\,G(2)}{G(6)\,H^{*}(1) - H(6)\,G^{*}(1)} = 
\frac{\eta(4\tau)\eta(6\tau)^{3}\eta(9\tau)\eta(24\tau)^{2}\eta(36\tau)\eta(144\tau)^{2}}{\eta(2\tau)\eta(3\tau)\eta(12\tau)^{2}\eta(18\tau)^{2}\eta(48\tau)^{2}\eta(72\tau)^{2}},\quad \Gamma_1(144), \quad -B=360.
\mylabel{eq:RR24a101}
\end{equation}

\subsubsection{$\chi_{24,3}$}
\mylabel{subsubsec:chi243}
We have 
\begin{align*}
G(1)&=G(1,24,\chi_{24,3})= \frac{1}
{
\eta_{24; 1,5}(\tau)
}
=
\frac{q^{-37/24}}
{
(q,q^{5},q^{19},q^{23};q^{24})_\infty
},\\
H(1)&=H(1,24,\chi_{24,3})= \frac{1}
{
\eta_{24; 7,11}(\tau)
}
=
\frac{q^{35/24}}
{
(q^{7},q^{11},q^{13},q^{17};q^{24})_\infty
}.
\end{align*}

\paragraph[0pt]{Type 1}
\mylabel{subsubsubsec:type24c.1}
\begin{alignat}{3}
G(2)\,H(1) - G(1)\,H(2) &= 
\frac{\eta(3\tau)\eta(4\tau)^{2}}{\eta(2\tau)\eta(24\tau)^{2}},&\quad &\Gamma_1(48), &\quad -B=24,
\mylabel{eq:RR24c11}\\
G(2)\,H(1) + G(1)\,H(2) &= 
\frac{\eta(6\tau)^{3}\eta(8\tau)\eta(12\tau)}{\eta(2\tau)\eta(3\tau)\eta(24\tau)^{3}},&\quad &\Gamma_1(48), &\quad -B=24,
\mylabel{eq:RR24c12}\\
G(3)\,H(1) - G(1)\,H(3) &= 
\frac{\eta(6\tau)^{2}\eta(8\tau)\eta(9\tau)\eta(36\tau)}{\eta(3\tau)\eta(18\tau)\eta(24\tau)^{2}\eta(72\tau)},&\quad &\Gamma_1(72), &\quad -B=72.
\mylabel{eq:RR24c13}
\end{alignat}

\paragraph[0pt]{Type 2}
\mylabel{subsubsubsec:type24c.2}
\begin{alignat}{3}
G(1)\,G(1) - H(1)\,H(1) &= 
\frac{\eta(2\tau)^{2}\eta(3\tau)^{2}\eta(8\tau)\eta(12\tau)^{3}}{\eta(\tau)^{2}\eta(6\tau)\eta(24\tau)^{5}},&\quad &\Gamma_1(24), &\quad -B=12.
\mylabel{eq:RR24c21}
\end{alignat}

\paragraph[0pt]{Type 6}
\mylabel{subsubsubsec:type24c.6}
\begin{alignat}{3}
G(1)\,H^{*}(1) - G^{*}(1)\,H(1) &= 
{2}\,\frac{\eta(6\tau)^{2}\eta(8\tau)^{2}\eta(12\tau)\eta(48\tau)^{3}}{\eta(2\tau)\eta(16\tau)\eta(24\tau)^{6}},&\quad &\Gamma_1(48), &\quad -B=24,
\mylabel{eq:RR24c61}\\
G(1)\,H^{*}(1) + G^{*}(1)\,H(1) &= 
{2}\,\frac{\eta(4\tau)\eta(8\tau)\eta(12\tau)^{3}\eta(48\tau)^{2}}{\eta(2\tau)\eta(24\tau)^{6}},&\quad &\Gamma_1(48), &\quad -B=24.
\mylabel{eq:RR24c62}
\end{alignat}

\paragraph[0pt]{Type 7}
\mylabel{subsubsubsec:type24c.7}
\begin{alignat}{3}
G^{*}(1)\,G(1) - H^{*}(1)\,H(1) = 
\frac{\eta(4\tau)\eta(6\tau)^{2}\eta(16\tau)\eta(24\tau)^{3}}{\eta(2\tau)\eta(8\tau)\eta(12\tau)^{2}\eta(48\tau)^{3}},\quad \Gamma_1(48), \quad -B=24,
\mylabel{eq:RR24c71}\\
G^{*}(1)\,G(1) + H^{*}(1)\,H(1) = 
\frac{\eta(4\tau)\eta(8\tau)\eta(12\tau)}{\eta(2\tau)\eta(48\tau)^{2}},\quad \Gamma_1(48), \quad -B=24.
\mylabel{eq:RR24c72}
\end{alignat}

\paragraph[0pt]{Type 8}
\mylabel{subsubsubsec:type24c.8}
\begin{alignat}{3}
G(1)^{2}\,H(2) - H(1)^{2}\,G(2) = 
{2}\,\frac{\eta(3\tau)\eta(6\tau)\eta(8\tau)^{2}\eta(12\tau)\eta(48\tau)^{3}}{\eta(\tau)\eta(16\tau)\eta(24\tau)^{6}},\quad \Gamma_1(48), \quad -B=24,
\mylabel{eq:RR24c81}\\
G(1)^{2}\,H(2) + H(1)^{2}\,G(2) = 
{2}\,\frac{\eta(3\tau)\eta(4\tau)\eta(8\tau)\eta(12\tau)^{3}\eta(48\tau)^{2}}{\eta(\tau)\eta(6\tau)\eta(24\tau)^{6}},\quad \Gamma_1(48), \quad -B=24.
\mylabel{eq:RR24c82}
\end{alignat}


\subsection{Mod $26$}
There is only one non-principal character mod $26$ that satisfies 
$\chi(-1)=1$,
namely the character $\chi_{26}$ induced by $\leg{\cdot}{13}$.
In this section
\begin{align*}
G(1)&=G(1,26,\chi_{26})= \frac{1}
{
\eta_{26; 1,3,9}(\tau)
}
=
\frac{q^{-7/4}}
{
(q,q^{3},q^{9},q^{17},q^{23},q^{25};q^{26})_\infty
},\\
H(1)&=H(1,26,\chi_{26})= \frac{1}
{
\eta_{26; 5,7,11}(\tau)
}
=
\frac{q^{5/4}}
{
(q^{5},q^{7},q^{11},q^{15},q^{19},q^{21};q^{26})_\infty
}.
\end{align*}
We find only one identity.

\subsubsection{Type 10}
\mylabel{subsubsec:type26.10}

\begin{equation}
\frac{G(3)\,H(2) - H(3)\,G(2)}{G(6)\,H^{*}(1) + H(6)\,G^{*}(1)} = 
\frac{\eta(26\tau)^{3}\eta(156\tau)^{3}}{\eta(52\tau)^{3}\eta(78\tau)^{3}},\quad \Gamma_1(156), \quad -B=576.
\mylabel{eq:RR26101}
\end{equation}

\subsection{Mod $28$}
\mylabel{subsec:mod28}
There is only one non-principal character mod $28$ that 
satisfies $\chi(-1)=1$,
namely the character $\chi(\cdot)=\leg{28}{\cdot}$.
In this section
\begin{align*}
G(1)&=G(1,28,\leg{28}{\cdot})= \frac{1}
{
\eta_{28; 1,3,9}(\tau)
}
=
\frac{q^{-17/8}}
{
(q,q^{3},q^{9},q^{19},q^{25},q^{27};q^{28})_\infty
},\\
H(1)&=H(1,28,\leg{28}{\cdot})= \frac{1}
{
\eta_{28; 5,11,13}(\tau)
}
=
\frac{q^{15/8}}
{
(q^{5},q^{11},q^{13},q^{15},q^{17},q^{23};q^{28})_\infty
}.
\end{align*}

\subsubsection{Type 1}
\mylabel{subsubsec:type28.1}
\begin{alignat}{3}
G(2)\,H(1) - G(1)\,H(2) &= 
\frac{\eta(4\tau)^{2}\eta(7\tau)\eta(14\tau)}{\eta(2\tau)\eta(28\tau)^{3}},&\quad &\Gamma_1(56), &\quad -B=48.
\mylabel{eq:RR2811}
\end{alignat}

\subsubsection{Type 6}
\mylabel{subsubsec:type28.6}
\begin{equation}
G(1)\,H^{*}(1) - G^{*}(1)\,H(1) = 
{2}\,\frac{\eta(4\tau)^{4}\eta(14\tau)^{3}\eta(56\tau)^{3}}{\eta(2\tau)^{2}\eta(8\tau)\eta(28\tau)^{7}},\quad \Gamma_1(56), \quad -B=48.
\mylabel{eq:RR2861}
\end{equation}

\subsubsection{Type 7}
\mylabel{subsubsec:type28.7}
\begin{equation}
G^{*}(1)\,G(1) - H^{*}(1)\,H(1) = 
\frac{\eta(4\tau)\eta(8\tau)\eta(28\tau)^{2}}{\eta(2\tau)\eta(56\tau)^{3}},\quad \Gamma_1(56), \quad -B=48.
\mylabel{eq:RR2871}
\end{equation}

\subsubsection{Type 8}
\mylabel{subsubsec:type28.8}
\begin{equation}
G(1)^{2}\,H(2) - H(1)^{2}\,G(2) = 
{2}\,\frac{\eta(4\tau)^{4}\eta(7\tau)\eta(14\tau)^{2}\eta(56\tau)^{3}}{\eta(\tau)\eta(2\tau)\eta(8\tau)\eta(28\tau)^{7}},\quad \Gamma_1(56), \quad -B=48.
\mylabel{eq:RR2881}
\end{equation}

\subsection{Mod $30$}
There is only one real non-principal character mod $30$ that satisfies 
$\chi(-1)=1$,
namely the character $\chi_{30}$ induced by the Legendre symbol mod $5$.
Thus in this section
\begin{align*}
G(1)&=G(1,30,\chi_{30})= \frac{1}
{
\eta_{30; 1,11}(\tau)
}
=
\frac{q^{-31/30}}
{
(q,q^{11},q^{19},q^{29};q^{30})_\infty
},\\
H(1)&=H(1,30,\chi_{30})= \frac{1}
{
\eta_{30; 7,13}(\tau)
}
=
\frac{q^{41/30}}
{
(q^{7},q^{13},q^{17},q^{23};q^{30})_\infty
}.
\end{align*}

\subsubsection{Type 6}
\mylabel{subsubsec:type30.6}
\begin{equation}
G(1)\,H^{*}(1) - G^{*}(1)\,H(1) = 
{2}\,\frac{\eta(4\tau)\eta(6\tau)^{2}\eta(60\tau)^{2}}{\eta(2\tau)\eta(12\tau)\eta(30\tau)^{3}},\quad \Gamma_1(60), \quad -B=48.
\mylabel{eq:RR3061}
\end{equation}

\subsubsection{Type 8}
\mylabel{subsubsec:type30.8}
\begin{equation}
G(1)^{2}\,H(2) - H(1)^{2}\,G(2) = 
{2}\,\frac{\eta(3\tau)\eta(4\tau)\eta(5\tau)\eta(6\tau)\eta(60\tau)^{2}}{\eta(\tau)\eta(10\tau)\eta(12\tau)\eta(15\tau)\eta(30\tau)^{2}},\quad \Gamma_1(60), \quad -B=48.
\mylabel{eq:RR3081}
\end{equation}

\subsection{Mod $34$}
There is only one real non-principal character mod $34$ that satisfies 
$\chi(-1)=1$,
namely the character $\chi_{34}$ induced by the Legendre symbol mod $17$.
\begin{align*}
G(1)&=G(1,34,\chi_{34})= \frac{1}
{
\eta_{34; 1,9,13,15}(\tau)
}
\\
&=
\frac{q^{2/3}}
{
(q,q^{9},q^{13},q^{15},q^{19},q^{21},q^{25},q^{33};q^{34})_\infty
},
\end{align*}
\begin{align*}
H(1)&=H(1,34,\chi_{34})= \frac{1}
{
\eta_{34; 3,5,7,11}(\tau)
}
\\
&=
\frac{q^{-4/3}}
{
(q^{3},q^{5},q^{7},q^{11},q^{23},q^{27},q^{29},q^{31};q^{34})_\infty
}.
\end{align*}

\subsubsection{Type 1}
\mylabel{subsubsec:type34.1}
\begin{alignat}{3}
G(2)\,H(1) - G(1)\,H(2) &= 
{-1}\,\frac{\eta(2\tau)^{2}\eta(17\tau)}{\eta(\tau)\eta(34\tau)^{2}},&\quad &\Gamma_1(68), &\quad -B=64.
\mylabel{eq:RR3411}
\end{alignat}

\subsubsection{Type 7}
\mylabel{subsubsec:type34.7}
\begin{equation}
G^{*}(1)\,G(1) - H^{*}(1)\,H(1) = 
{-1}\,\frac{\eta(4\tau)}{\eta(68\tau)},\quad \Gamma_1(68), \quad -B=64.
\mylabel{eq:RR3471}
\end{equation}

\subsubsection{Type 9}
\mylabel{subsubsec:type34.9}
\begin{equation}
G(1)^{2}\,H(1) - H(1)^{2}\,G(1) = 
{-1}\,\frac{\eta(2\tau)^{2}\eta(17\tau)}{\eta(\tau)\eta(34\tau)^{2}},\quad \Gamma_1(34), \quad -B=16.
\mylabel{eq:RR3491}
\end{equation}

\subsection{Mod $40$}
There are three real non-principal characters mod $40$ that satisfy   
$\chi(-1)=1$.

\begin{enumerate}
\item[(i)]
The character $\chi_{40,1}(\cdot)$ induced by $\leg{\cdot}{5}$.
This actually a character mod $10$. See Section \subsect{mod10}.
\item[(ii)]
The character $\chi_{40,2}(\cdot)$ induced by $\leg{8}{\cdot}$.
\item[(iii)]
The character $\chi_{40,3}(\cdot)=\leg{40}{\cdot}$.
\end{enumerate}

\subsubsection{$\chi_{40,2}$}
\mylabel{subsubsec:chi402}
\begin{align*}
G(1)&=G(1,40,\chi_{40,2})= \frac{1}
{
\eta_{40; 1,7,9,17}(\tau)
}
\\
&=
\frac{q^{-19/12}}
{
(q,q^{7},q^{9},q^{17},q^{23},q^{31},q^{33},q^{39};q^{40})_\infty
},
\end{align*}
\begin{align*}
H(1)&=H(1,40,\chi_{40,2})= \frac{1}
{
\eta_{40; 3,11,13,19}(\tau)
}
\\
&=
\frac{q^{17/12}}
{
(q^{3},q^{11},q^{13},q^{19},q^{21},q^{27},q^{29},q^{37};q^{40})_\infty
}.
\end{align*}

\paragraph[0pt]{Type 1}
\mylabel{subsubsubsec:type40b.1}
\begin{alignat}{3}
G(2)\,H(1) - G(1)\,H(2) &= 
\frac{\eta(4\tau)^{2}\eta(10\tau)^{2}}{\eta(2\tau)\eta(8\tau)\eta(20\tau)\eta(40\tau)},&\quad &\Gamma_1(80), &\quad -B=80.
\mylabel{eq:RR40b11}
\end{alignat}

\paragraph[0pt]{Type 6}
\mylabel{subsubsubsec:type40b.6}
\begin{equation}
G(1)\,H^{*}(1) + G^{*}(1)\,H(1) = 
{2}\,\frac{\eta(8\tau)^{3}\eta(10\tau)\eta(20\tau)^{3}\eta(80\tau)^{3}}{\eta(2\tau)\eta(16\tau)\eta(40\tau)^{8}},\quad \Gamma_1(80), \quad -B=112.
\mylabel{eq:RR40b61}
\end{equation}

\paragraph[0pt]{Type 7}
\mylabel{subsubsubsec:type40b.7}
\begin{equation}
G^{*}(1)\,G(1) + H^{*}(1)\,H(1) = 
\frac{\eta(4\tau)\eta(8\tau)\eta(10\tau)}{\eta(2\tau)\eta(16\tau)\eta(80\tau)},\quad \Gamma_1(80), \quad -B=80.
\mylabel{eq:RR40b71}
\end{equation}

\paragraph[0pt]{Type 8}
\mylabel{subsubsubsec:type40b.8}
\begin{equation}
G(1)^{2}\,H(2) + H(1)^{2}\,G(2) = 
{2}\,\frac{\eta(4\tau)^{2}\eta(5\tau)\eta(16\tau)\eta(20\tau)\eta(80\tau)}{\eta(\tau)\eta(8\tau)^{2}\eta(40\tau)^{3}},\quad \Gamma_1(80), \quad -B=80.
\mylabel{eq:RR40b81}
\end{equation}

\subsubsection{$\chi_{40,3}$}
\mylabel{subsubsec:chi403}
\begin{align*}
G(1)&=G(1,40,\chi_{40,3})= \frac{1}
{
\eta_{40; 1,3,9,13}(\tau)
}
\\
&=
\frac{q^{-43/12}}
{
(q,q^{3},q^{9},q^{13},q^{27},q^{31},q^{37},q^{39};q^{40})_\infty
},
\end{align*}
\begin{align*}
H(1)&=H(1,40,\chi_{40,3})= \frac{1}
{
\eta_{40; 7,11,17,19}(\tau)
}
\\
&=
\frac{q^{41/12}}
{
(q^{7},q^{11},q^{17},q^{19},q^{21},q^{23},q^{29},q^{33};q^{40})_\infty
}.
\end{align*}

\paragraph[0pt]{Type 1}
\mylabel{subsubsubsec:type40c.1}
\begin{alignat}{3}
G(2)\,H(1) - G(1)\,H(2) &= 
\frac{\eta(4\tau)\eta(8\tau)\eta(10\tau)^{2}}{\eta(2\tau)\eta(40\tau)^{3}},&\quad &\Gamma_1(80), &\quad -B=112.
\mylabel{eq:RR40c11}
\end{alignat}

\paragraph[0pt]{Type 6}
\mylabel{subsubsubsec:type40c.6}
\begin{equation}
G(1)\,H^{*}(1) + G^{*}(1)\,H(1) = 
{2}\,\frac{\eta(8\tau)^{3}\eta(10\tau)\eta(20\tau)^{3}\eta(80\tau)^{3}}{\eta(2\tau)\eta(16\tau)\eta(40\tau)^{8}},\quad \Gamma_1(80), \quad -B=112.
\mylabel{eq:RR40c61}
\end{equation}

\paragraph[0pt]{Type 7}
\mylabel{subsubsubsec:type40c.7}
\begin{equation}
G^{*}(1)\,G(1) + H^{*}(1)\,H(1) = 
\frac{\eta(4\tau)\eta(10\tau)\eta(16\tau)\eta(40\tau)}{\eta(2\tau)\eta(80\tau)^{3}},\quad \Gamma_1(80), \quad -B=112.
\mylabel{eq:RR40c71}
\end{equation}

\paragraph[0pt]{Type 8}
\mylabel{subsubsubsec:type40c.8}
\begin{equation}
G(1)^{2}\,H(2) + H(1)^{2}\,G(2) = 
{2}\,\frac{\eta(5\tau)\eta(8\tau)^{3}\eta(20\tau)^{3}\eta(80\tau)^{3}}{\eta(\tau)\eta(16\tau)\eta(40\tau)^{8}},\quad \Gamma_1(80), \quad -B=112.
\mylabel{eq:RR40c81}
\end{equation}

\subsection{Mod $42$}
There is only one real non-principal character mod $42$ that satisfies 
$\chi(-1)=1$,
namely the one induced by the mod $21$ character 
$\chi_{42}(\cdot) = \leg{\cdot}{3}\,\leg{\cdot}{7}$.

In this section
\begin{align*}
G(1)&=G(1,42,\chi_{42})= \frac{1}
{
\eta_{42; 1,5,17}(\tau)
}
=
\frac{q^{-11/4}}
{
(q,q^{5},q^{17},q^{25},q^{37},q^{41};q^{42})_\infty
},\\
H(1)&=H(1,42,\chi_{42})= \frac{1}
{
\eta_{42; 11,13,19}(\tau)
}
=
\frac{q^{13/4}}
{
(q^{11},q^{13},q^{19},q^{23},q^{29},q^{31};q^{42})_\infty
}.
\end{align*}

\subsubsection{Type 1}
\mylabel{subsubsec:type42.1}
\begin{alignat}{3}
G(2)\,H(1) - G(1)\,H(2) &= 
\frac{\eta(4\tau)\eta(6\tau)^{2}\eta(7\tau)}{\eta(2\tau)\eta(12\tau)\eta(21\tau)\eta(42\tau)},&\quad &\Gamma_1(84), &\quad -B=96.
\mylabel{eq:RR4211}
\end{alignat}

\subsubsection{Type 7}
\mylabel{subsubsec:type42.7}
\begin{equation}
G^{*}(1)\,G(1) - H^{*}(1)\,H(1) = 
\frac{\eta(4\tau)\eta(6\tau)\eta(28\tau)}{\eta(2\tau)\eta(84\tau)^{2}},\quad \Gamma_1(84), \quad -B=96.
\mylabel{eq:RR4271}
\end{equation}

\subsection{Mod $56$}
There are three real non-principal characters mod $56$ that satisfy   
$\chi(-1)=1$.
\begin{enumerate}
\item[(i)]
The character $\leg{56}{\cdot}$.                                            
\item[(ii)]
The character induced by the mod $28$ character $\leg{28}{\cdot}$.  
See Section \subsect{mod28}.
\item[(iii)]
The character induced by the mod $8$ character $\leg{8}{\cdot}$.  
\end{enumerate}
Only the third character led to new identities. In this section
we assume $\chi$ is the mod $56$ character induced by $\leg{8}{\cdot}$.
Thus in this section
\begin{align*}
G(1,56,\chi)&=G(1) = \frac{1}
{
\eta_{56; 1,9,15,17,23,25}(\tau)
}
\\
&=
\frac{q^{11/8}}{
(q,q^{9},q^{15},q^{17},q^{23},q^{25},q^{31},q^{33},q^{39},q^{41},q^{47},q^{55};q^{56})_\infty
},
\end{align*}
\begin{align*}
H(1,56,\chi)&=H(1) = \frac{1}
{
\eta_{56; 3,5,11,13,19,27}(\tau)
}
\\
&= 
\frac{q^{-13/8}}{                
(q^{3},q^{5},q^{11},q^{13},q^{19},q^{27},q^{29},q^{37},q^{43},q^{45},q^{51},q^{53};q^{56})_\infty
}
\end{align*}

\subsubsection{Type 1}
\mylabel{subsubsec:type56c.1}
\begin{alignat}{3}
G(2)\,H(1) + G(1)\,H(2) &= 
\frac{\eta(2\tau)\eta(4\tau)\eta(14\tau)}{\eta(\tau)\eta(8\tau)\eta(56\tau)},&\quad &\Gamma_1(112), &\quad -B=144.
\mylabel{eq:RR56c11}
\end{alignat}

\subsubsection{Type 6}
\mylabel{subsubsec:type56c.6}
\begin{equation}
G(1)\,H^{*}(1) - G^{*}(1)\,H(1) = 
{2}\,\frac{\eta(4\tau)\eta(14\tau)\eta(16\tau)\eta(112\tau)}{\eta(8\tau)^{2}\eta(56\tau)^{2}},\quad \Gamma_1(112), \quad -B=144.
\mylabel{eq:RR56c61}
\end{equation}

\subsubsection{Type 7}
\mylabel{subsubsec:type56c.7}
\begin{equation}
G^{*}(1)\,G(1) - H^{*}(1)\,H(1) = 
{-1}\,\frac{\eta(8\tau)\eta(14\tau)\eta(56\tau)}{\eta(16\tau)\eta(28\tau)\eta(112\tau)},\quad \Gamma_1(112), \quad -B=144.
\mylabel{eq:RR56c71}
\end{equation}

\subsubsection{Type 8}
\mylabel{subsubsec:type56c.8}
\begin{equation}
G(1)^{2}\,H(2) - H(1)^{2}\,G(2) = 
{2}\,\frac{\eta(2\tau)\eta(4\tau)\eta(7\tau)\eta(16\tau)\eta(112\tau)}{\eta(\tau)\eta(8\tau)^{2}\eta(56\tau)^{2}},\quad \Gamma_1(112), \quad -B=144.
\mylabel{eq:RR56c81}
\end{equation}

\subsection{Mod $60$}
There are three real non-principal characters mod $60$ that satisfy   
$\chi(-1)=1$.
\begin{enumerate}
\item[(i)]
The character induced by $\leg{\cdot}{5}$.
\item[(ii)]
The character $\chi_{60,2}(\cdot) = \leg{60}{\cdot}$.                                
See Section ??.
\item[(iii)]
The character $\chi_{60,3}(\cdot)$ induced by the mod $12$ character 
$\leg{12}{\cdot}$.  
\end{enumerate}
Only (ii), (iii) seem to lead to new identities.

\subsubsection{$\chi_{60,2}$}
\mylabel{subsubsec:chi602}
In this section
\begin{align*}
G(1,60, \chi_{60,2}) &= G(1) =
\frac{1}{
\eta_{60; 1,7,11,17}(\tau)
}
\\
&=
\frac{q^{-35/6}}
{
(q,q^{7},q^{11},q^{17},q^{43},q^{49},q^{53},q^{59};q^{60})_\infty
},
\end{align*}
\begin{align*}
H(1,60, \chi_{60,2}) &= H(1) =
\frac{1}{
\eta_{60; 13,19,23,29}(\tau)
}
\\
&=
\frac{q^{37/6}}
{
(q^{13},q^{19},q^{23},q^{29},q^{31},q^{37},q^{41},q^{47};q^{60})_\infty
}.
\end{align*}

\paragraph[0pt]{Type 1}
\mylabel{subsubsubsec:type60b.1}
\begin{alignat}{3}
G(2)\,H(1) - G(1)\,H(2) &= 
\frac{\eta(4\tau)\eta(6\tau)\eta(10\tau)\eta(30\tau)}{\eta(2\tau)\eta(60\tau)^{3}},&\quad &\Gamma_1(120), &\quad -B=192.
\mylabel{eq:RR60b11}
\end{alignat}

\paragraph[0pt]{Type 6}
\mylabel{subsubsubsec:type60b.6}
\begin{equation}
G(1)\,H^{*}(1) - G^{*}(1)\,H(1) = 
{2}\,\frac{\eta(4\tau)\eta(12\tau)^{2}\eta(20\tau)^{2}\eta(30\tau)^{3}\eta(120\tau)^{4}}{\eta(2\tau)\eta(24\tau)\eta(40\tau)\eta(60\tau)^{9}},\quad \Gamma_1(120), \quad -B=192.
\mylabel{eq:RR60b61}
\end{equation}

\paragraph[0pt]{Type 7}
\mylabel{subsubsubsec:type60b.7}
\begin{equation}
G^{*}(1)\,G(1) - H^{*}(1)\,H(1) = 
\frac{\eta(4\tau)\eta(6\tau)\eta(10\tau)\eta(24\tau)\eta(40\tau)\eta(60\tau)^{3}}{\eta(2\tau)\eta(12\tau)\eta(20\tau)\eta(30\tau)\eta(120\tau)^{4}},\quad \Gamma_1(120), \quad -B=192.
\mylabel{eq:RR60b71}
\end{equation}

\paragraph[0pt]{Type 8}
\mylabel{subsubsubsec:type60b.8}
\begin{equation}
G(1)^{2}\,H(2) - H(1)^{2}\,G(2) = 
{2}\,\frac{\eta(3\tau)\eta(4\tau)\eta(5\tau)\eta(12\tau)^{2}\eta(20\tau)^{2}\eta(30\tau)^{4}\eta(120\tau)^{4}}{\eta(\tau)\eta(6\tau)\eta(10\tau)\eta(15\tau)\eta(24\tau)\eta(40\tau)\eta(60\tau)^{9}},\quad \Gamma_1(120), \quad -B=192.
\mylabel{eq:RR60b81}
\end{equation}

\subsubsection{$\chi_{60,3}$}
\mylabel{subsubsec:chi603}
In this section
\begin{align*}
G(1,60, \chi_{60,3}) &= G(1) =
\frac{1}{
\eta_{60; 1,11,13,23}(\tau)
}
\\
&=
\frac{q^{-17/6}}
{
(q,q^{11},q^{13},q^{23},q^{37},q^{47},q^{49},q^{59};q^{60})_\infty
},
\end{align*}
\begin{align*}
H(1,60, \chi_{60,3}) &= H(1) =
\frac{1}{
\eta_{60; 7,17,19,29}(\tau)
}
\\
&=\frac{q^{19/6}}{
(q^{7},q^{17},q^{19},q^{29},q^{31},q^{41},q^{43},q^{53};q^{60})_\infty
}.
\end{align*}

\paragraph[0pt]{Type 1}
\mylabel{subsubsubsec:type60c.1}
\begin{alignat}{3}
G(2)\,H(1) - G(1)\,H(2) &= 
\frac{\eta(4\tau)\eta(6\tau)^{2}\eta(10\tau)}{\eta(2\tau)\eta(12\tau)^{2}\eta(60\tau)},&\quad &\Gamma_1(120), &\quad -B=160.
\mylabel{eq:RR60c11}
\end{alignat}

\paragraph[0pt]{Type 6}
\mylabel{subsubsubsec:type60c.6}
\begin{equation}
G(1)\,H^{*}(1) - G^{*}(1)\,H(1) = 
{2}\,\frac{\eta(4\tau)\eta(6\tau)^{2}\eta(20\tau)^{2}\eta(24\tau)\eta(30\tau)\eta(120\tau)^{2}}{\eta(2\tau)\eta(12\tau)^{3}\eta(40\tau)\eta(60\tau)^{4}},\quad \Gamma_1(120), \quad -B=160.
\mylabel{eq:RR60c61}
\end{equation}

\paragraph[0pt]{Type 7}
\mylabel{subsubsubsec:type60c.7}
\begin{equation}
G^{*}(1)\,G(1) - H^{*}(1)\,H(1) = 
\frac{\eta(4\tau)\eta(6\tau)\eta(10\tau)\eta(40\tau)\eta(60\tau)^{2}}{\eta(2\tau)\eta(20\tau)\eta(24\tau)\eta(30\tau)\eta(120\tau)^{2}},\quad \Gamma_1(120), \quad -B=160.
\mylabel{eq:RR60c71}
\end{equation}

\paragraph[0pt]{Type 8}
\mylabel{subsubsubsec:type60c.8}
\begin{equation}
G(1)^{2}\,H(2) - H(1)^{2}\,G(2) = 
{2}\,\frac{\eta(3\tau)\eta(4\tau)\eta(5\tau)\eta(6\tau)\eta(20\tau)^{2}\eta(24\tau)\eta(30\tau)^{2}\eta(120\tau)^{2}}{\eta(\tau)\eta(10\tau)\eta(12\tau)^{3}\eta(15\tau)\eta(40\tau)\eta(60\tau)^{4}},\quad \Gamma_1(120), \quad -B=160.
\mylabel{eq:RR60c81}
\end{equation}


\providecommand{\bysame}{\leavevmode\hbox to3em{\hrulefill}\thinspace}
\providecommand{\MR}{\relax\ifhmode\unskip\space\fi MR }
\providecommand{\MRhref}[2]{%
  \href{http://www.ams.org/mathscinet-getitem?mr=#1}{#2}
}
\providecommand{\href}[2]{#2}

\end{document}